\documentclass[final,12pt]{elsarticle}
\usepackage{graphicx}
\usepackage{epstopdf}
\usepackage{amsmath}
\usepackage{amssymb}
\usepackage{enumerate}
\usepackage{enumitem}
\usepackage{float}

\newcommand{\req}[1]{Eq.\,(\ref{#1})}

\numberwithin{equation}{section}

\begin{document} 
\begin{frontmatter}
\title{Boltzmann Equation Solver Adapted to Emergent Chemical Non-equilibrium}
\author{Jeremiah Birrell$^{a,b}$, Jon Wilkening$^{c}$, Johann Rafelski$^{b}$,\\
\scriptsize{$^{a}$Program in Applied Mathematics, and $^{b}$Department of Physics\\ The University of Arizona, Tucson, Arizona, 85721, USA\\ $^{c}$ Department of Mathematics and Lawrence Berkeley National Laboratory,\\
 University of California, Berkeley, CA 94721 }\\
}

\begin{abstract}
We present a novel method to solve the spatially homogeneous and isotropic relativistic Boltzmann equation. We employ a basis set of orthogonal polynomials dynamically adapted to allow for emergence of chemical non-equilibrium. Two time dependent parameters characterize the set of orthogonal polynomials, the effective temperature $T(t)$ and  phase space occupation factor $\Upsilon(t)$. In this first paper we address  (effectively) massless fermions and derive dynamical equations for $T(t)$ and $\Upsilon(t)$ such that the zeroth order term of the basis alone captures the particle number density and energy density of each particle distribution. We validate our method and illustrate the reduced computational cost and the ability to easily represent final state chemical non-equilibrium by studying a model problem that is motivated by the physics of the  neutrino freeze-out process in the early Universe, where  the essential physical characteristics  include reheating from another disappearing particle component ($e^\pm$-annihilation). 
\end{abstract}

\end{frontmatter}
\section{Introduction}
\subsection{Different Equilibria}
At sufficiently high temperatures, such as existed in the early universe and reproduced in modern day relativistic heavy ion collider studies of  the hot, ultra-dense state of matter called quark gluon plasma (QGP), both particle creation and annihilation (i.e. chemical) processes and momentum exchanging (i.e. kinetic) scattering processes can occur sufficiently rapidly to establish complete thermal equilibrium. The most probable canonical distribution function $f_{ch}^\pm$ of  fermions (+) and bosons (-) in both chemical and kinetic equilibrium is found by maximizing microcanonical entropy subject to energy being conserved
\begin{equation}\label{ch_eq}
f_{ch}^\pm=\frac{1}{\exp(E/T)\pm 1}, \hspace{2mm} T>T_{ch}
\end{equation}
where $E$ is the particle energy, $T$ the temperature, and $T_{ch}$ the chemical freeze-out temperature.

For a physical system comprising {\em interacting} particles whose temperature is decreasing with time, there will be a period where the temperature is greater than the kinetic freeze-out temperature, $T_k$, but below chemical freeze-out. During this period, momentum exchanging processes continue to maintain an equilibrium distribution of energy among the available particles, which we call kinetic equilibrium, but particle number changing processes no longer occur rapidly enough to keep the equilibrium particle number yield.  For $T<T_{ch}$ the particle number changing processes have `frozen-out'. In this condition the momentum distribution, which is in kinetic equilibrium but chemical non-equilibrium, is obtained by maximizing  microcanonical entropy subject to  particle number and energy constraints and thus two parameters appear
\begin{equation}\label{k_eq}
f_{k}^\pm=\frac{1}{\Upsilon^{-1} \exp(E/T)\pm 1},\hspace{2mm} T_k<T\leq T_{ch}.
\end{equation}
The need to preserve the total particle number within the distribution introduces an additional parameter $\Upsilon$ called fugacity. 

A fugacity different from $1$ implies an over-abundance ($\Upsilon>1$) or under-abundance ($\Upsilon<1$) of particles compared to chemical equilibrium and in either of these  situations one speaks of chemical non-equilibrium. We emphasize that in our context, the appearance of a dynamical fugacity $\Upsilon\ne 1$ is  not related to conserved quantum numbers. Conserved quantum numbers such as baryon number and charge introduce, through chemical potentials, fugacities that enhance e.g particles and deplete antiparticles. In the presence of conservation laws one has for particles $\Upsilon_p=\gamma e^{\mu_p/T}$ and antiparticles  $\Upsilon_a=\gamma e^{-\mu_a/T}$ where the quantity $\gamma$ is the dynamical fugacity. For fermions (+) the distribution \req{k_eq} always satisfies the Fermi-Dirac constraint $f_{k}^+\le 1$, however the distribution of particles can approach much more closely the quantum degeneracy condition  $f_{k}^+\to 1$ for $\Upsilon \gg  1$.

Once the temperature drops below the kinetic freeze-out temperature $T_k$ we reach  the free streaming period where  particle scattering processes have completely frozen out and the resultant distribution is obtained by solving the collisionless Boltzmann equation with initial condition as given by the chemical non-equilibrium   distribution \req{k_eq}.  As already indicated, the two transitions between these three regimes constitute  the freeze-out process -- first we have at $T_{ch}$ the chemical freeze-out and at lower $T_k$ the kinetic freeze-out.

\subsection{Boltzmann Evolution and Chemical Non-Equilibrium}
Exact chemical and kinetic equilibrium and sharp freeze-out transitions at $T_{ch}$ and $T_k$ are  only approximations.  The  Boltzmann equation is a more precise model of the dynamics of the freeze-out process and furthermore, given the collision dynamics it is capable of capturing in a {\em quantitative manner} the non-thermal distortions from equilibrium, for example the emergence of actual distributions and the approximate values  of $T_{ch}$, $T_k$, and $\Upsilon$.  Indeed,  in  such a dynamical description no hypothesis about the presence of kinetic or chemical (non) equilibrium needs to be made, as the distribution similar to \req{k_eq} with   $\Upsilon\ne  1$ emerges naturally as the outcome of collision processes, even when the particle system approaches the freeze-out temperature domain  in chemical equilibrium.

Considering this physical situation it is striking that the literature on Boltzmann solvers does not reflect on the accommodation of emergent chemical non-equilibrium into the method of solution. For an all-numerical solver this may not be a necessary step as long as there are no constraints that preclude development of a general non-equilibrium solution. However, when strong chemical non-equilibrium is present either in the intermediate time period or/and at the end of the evolution a brute force approach can be very costly in computer time. Motivated by this circumstance and past work with physical environments in which chemical non-equilibrium arose,  we introduce here a  spectral method for solving the Boltzmann equation that utilizes a dynamical basis of orthogonal polynomials which is adapted to the case of emerging chemical non-equilibrium. We validate our method via a  model problem  that captures the essential physical characteristics of interest and use it to highlight the type of situation where this new method exhibits its advantages.

In the cosmological neutrino freeze-out context, the general relativistic Boltzmann equation has been used to study neutrino freeze-out in the early universe and has been successfully solved using both discretization in momentum space \cite{Madsen,Dolgov_Hansen,Gnedin,Mangano2005} and a spectral method based on a fixed basis of orthogonal polynomials \cite{Esposito2000,Mangano2002}.    In Refs.\cite{Wilkening,Wilkening2} the non-relativistic Boltzmann equation was solved via a spectral method similar in  one important mathematical idea to the approach we present here.  For near equilibrium solutions, the spectral methods have the advantage of requiring a relatively small number of modes to obtain an accurate solution, as opposed to momentum space discretization which in general leads to a large highly coupled nonlinear system of odes irrespective of the near equilibrium nature of the system.  

The efficacy of the spectral method used in \cite{Esposito2000,Mangano2002} can largely be attributed to the fact that, under the conditions considered there, the true solution is very close to a chemical equilibrium distribution \req{ch_eq} where the temperature is controlled by the dilution of the system. However, the recent PLANCK CMB results \cite{Planck} indicate the possibility that neutrinos participated in reheating to a greater degree than previously believed. As we discussed recently~\cite{Birrell} this can also lead to a more pronounced chemical non-equilibrium. Efficiently obtaining this emergent chemical non-equilibrium within realm of kinetic theory motivates the development of a new numerical method that adapts to this new circumstance.

With this novel Boltzmann solver we present here it is also possible to return to the exploration of the emergent chemical non-equilibrium in processes governing  laboratory QGP physics. Indeed,  the study of chemical non-equilibrium as a separate process from kinetic equilibrium has its roots in the field of laboratory  QGP formation and observation, of which one of the proposed observables is the newly produced strange quark flavor~\cite{Muller:2011tu}.  The  chemical non-equilibrium  analysis method~\cite{Rafelski:1991} is today the only successful  statistical hadronization model for experimental results~\cite{Petran:2013dva} confirming chemical non-equilibrium for all  strongly interacting particles produced by a QGP fireball~\cite{Petran:2013lja}. 

This paper establishes our new method for the case of ultra-relativistic particles does not address  the pertinent physical applications in neutrino cosmology or quark-gluon plasma physics.  The former is treated in our paper \cite{Birrell_nu_param} and the later will be a subject of future work. While in this  work we address the case  where the mass-scale  of particles is entirely irrelevant, we have also developed a similar method for the case that the particle mass is non-negligible. Introduction of a mass scale presents no major conceptual modifications, but requires detailed technical modifications from the simpler scheme we present here  that don't lend themselves well to a simultaneous presentation of both methods.  The method including mass will be introduced and addressed in the context of specific applications of the method in future publications.

In section \ref{boltzmann_basics} we give a basic overview of the relativistic Boltzmann equation in the format aiming to address the early Universe neutrino freeze-out process, but which can be easily recast into the format appropriate for other applications.  In section \ref{the_method} we discuss our modified spectral method in detail.  In subsection \ref{free_stream_approach} we recall the orthogonal polynomial basis used in \cite{Esposito2000,Mangano2002} and in subsection \ref{kinetic_eq_approach} we introduce our modified basis and characterize precisely the differences in the method  we propose. We compare these two bases in subsection \ref{basis_comparison}. In subsection \ref{dynamics_sec} we use the Boltzmann equation to derive the dynamics of the mode coefficients and identify physically motivated evolution equations for the effective temperature and fugacity.  In section \ref{validation} we validate the method using a model problem.  In  \ref{orthopoly_app} we give further details on the construction of the parametrized family orthogonal polynomials we use to solve the Boltzmann equation.  

\section{Relativistic Boltzmann Equation }\label{boltzmann_basics}
Consider the relativistic Bolzmann equation for systems of fermions under the assumption of homogeneity and isotropy. We assume that the particle are effectively massless  i.e. the temperature is much greater than the mass scale and use the sign convention $(+,-,-,-)$.   Homogeneity and isotropy in a given reference frame imply that the distribution function of each particle species under consideration has the form $f^+=f^+(t,p)$ in that frame. Here $p$ is  the magnitude of the relativistic momentum, related to the relativistic energy $E$ by the usual relation $E=\sqrt{p^2+m^2}$, where $m$ is the particle mass.  Under such assumptions, the Boltzmann equation reduces to
\begin{equation}\label{boltzmann_p}
\partial_t f^+-pH \partial_p f^+=\frac{1}{E}C[f^+],\hspace{2mm} H=\frac{\dot{a}}{a}
\end{equation}
where $f^+$ is the particle distribution function and we drop from now on the upper Fermi index (+) .   

In \req{boltzmann_p} we have allowed for a dilution of the system, encoded in the linear scale factor $a(t)$.  In cosmological applications it represents the scale factor of the universe and $H$ a dilution rate called the Hubble parameter.  We do not specify its dynamics here, as that is application specific.  When validating our method, we will use an externally prescribed function motivated by applications to cosmology.  For a more in depth discussion of the relativistic Boltzmann equation see for example \cite{Cercignani,Ehlers}. 

The term $C[f]$ on the right hand side of the Boltzmann equation is called the collision operator and models the short range scattering processes that cause deviations from geodesic motion. For $2\leftrightarrow 2$ reactions between fermions, such as neutrinos, the collision operator takes the form
\begin{align}\label{coll}
C[f_1]=&\frac{1}{2}\int F(p_1,p_2,p_3,p_4) S |\mathcal{M}|^2(2\pi)^4\delta(\Delta p)\prod_{i=2}^4\frac{d^{3}p_i}{2(2\pi)^3E_i},\\
F=&f_3(p_3)f_4(p_4)f^1(p_1)f^2(p_2)-f_1(p_1)f_2(p_2)f^3(p_3)f^4(p_4),\hspace{2mm} f^i=1- f_i.\notag
\end{align}
Here $|\mathcal{M}|^2$ is the process amplitude and $S$ is a numerical factor that incorporates symmetries and prevents over-counting.  There are many references that discuss the matrix elements required to simulate various application domains, such as neutrino freeze-out \cite{Madsen,Dolgov_Hansen}, or strong interaction processes  in QGP and hadron gas~\cite{Letessier:2002gp,Kuznetsova:2010pi}.  In the analysis of our method in this work, we will use an idealized version of the collision operator in order to avoid the details of any particular application domain.  For an example of an application of the here presented method to a physically realistic problem see \cite{Birrell_nu_param}, where we study the freeze-out of cosmological neutrinos.

The distribution $f$  is normalized so that the particle number density is given by
\begin{equation}\label{n+dens}
n=\frac{g_p}{(2\pi)^3}\int  f(p) d^3p.
\end{equation}
where $g_p$ is the degeneracy of the particle species. The energy density and pressure associated with the distribution function $f$ are obtained from the stress energy tensor
\begin{equation}
T^{\mu\nu}=\frac{g_p}{(2\pi)^3}\int \frac{p^\mu p^\nu}{E} f(p) d^3p.
\end{equation}
 We emphasize that all momenta appearing in this work are the physical momenta and not, for example, the coordinate momenta in coordinates where the cosmological metric tensor takes the form $ds^2=dt^2-a(t)^2(dx^2+dy^2+dz^2)$ (i.e. we work in an orthonormal basis for the tangent spaces). In particular this means the usual relation familiar from special relativity, $p^0=E=\sqrt{m^2+p^2}$, holds.  The only essential difference between using these results in a special relativistic versus cosmological context is whether or not $a(t)$ is allowed to change in time. Because of this, the energy density $\rho=T^{00}$ is
\begin{equation}\label{rho_dens}
\rho=\frac{g_p}{(2\pi)^3}\int E f(p) d^3p
\end{equation}

The Boltzmann equation \req{boltzmann_p} can be simplified by the method of characteristics. Writing $f(p; t)=g(a(t)p,t)$ and reverting back to call the new distribution $g\to f$, the 2nd term in \req{boltzmann_p} cancels out and the evolution in time can be studied directly.  This transformation implies for the rate of change in the  particle number density and energy density  
\begin{align}\label{n_div}
\frac{1}{a^3}\frac{d}{dt}(a^3n_1)=&\frac{g_p}{(2\pi)^3}\int C[f_1] \frac{d^3p}{E}.\\
\label{rho_div}
\frac{1}{a^4}\frac{d}{dt}(a^4\rho_1)=&\frac{g_p}{(2\pi)^3}\int C[f_1] d^3p .
\end{align} 
For free-streaming particles the vanishing of the collision operator implies conservation of `comoving' particle number of species 1. From the associated powers of $a$ in \req{n_div} and \req{rho_div} we see that the  energy per free streaming particle as measured by an observer scales as $1/a$. 

\section{Spectral Methods}\label{the_method}
\subsection{Polynomials for systems close to kinetic and chemical equilibrium}\label{free_stream_approach}
Here we outline the approach for solving \req{a_vars} used in \cite{Esposito2000,Mangano2002} in order to contrast it with our approach as presented in subsection \ref{kinetic_eq_approach}.  As just discussed, the Boltzmann equation  is a linear first order partial differential equation and can be reduced using a new variable $y=a(t)p$  via the method of characteristics and exactly solved in the collision free ($C[f]=0)$ limit.   This motivates a change of variables from $p$ to $y$ which eliminates the momentum derivative, leaving the simplified equation
\begin{equation}\label{a_vars}
\partial_tf=\frac{1}{E} C[f].
\end{equation}

We let $\hat\chi_i$ be the orthonormal polynomial basis on the interval $[0,\infty)$ with respect to the weight function
\begin{equation}\label{free_stream_weight}
f_{ch}=\frac{1}{e^y+1},
\end{equation}
constructed as in  \ref{orthopoly_app}. $f_{ch}$ is the Fermi-Dirac chemical equilibrium distribution for massless fermions and temperature $T=1/a$.  Therefore this ansatz is well suited to distributions that are manifestly in chemical equilibrium ($\Upsilon=1$) or remain close and with $T\propto 1/a$, which we call dilution temperature scaling.  Assuming that $f$ is such a distribution  motivates the decomposition
\begin{equation}\label{free_stream_ansatz}
f=f_{ch}\chi,\qquad \chi=\sum_i d^i\hat\chi_i.
\end{equation}

Using this ansatz  equation \req{a_vars} becomes
\begin{equation}
\dot{d}^k=\int_0^\infty\frac{1}{E}\hat{\chi}_k C[f]dy.
\end{equation}
Because of \req{free_stream_ansatz}, we call this the chemical equilibrium method.

We also have the following expressions for the particle number density and energy density
\begin{align}\label{free_stream_moments}
n&=\frac{g_p}{2\pi^2 a^3}\sum_0^2 d^i\int_0^\infty f_{ch}\hat\chi_i y^2dy,\\
\rho&=\frac{g_p}{2\pi^2a^4}\sum_0^3 d^i\int_0^\infty f_{ch}\hat\chi_i y^3dy.
\end{align}

Note that the sums truncate at $3$ and $4$ terms respectively, due to the fact that $\hat\chi_k$ is orthogonal to all polynomials of degree less than $k$. This implies that in general, at least four modes are required to capture both the particle number and energy flow. More modes are needed if the non-thermal distortions are large and the back reaction of higher modes on lower modes is significant.

\subsection{Polynomials for systems   not close to chemical equilibrium}\label{kinetic_eq_approach}
Our primary interest is in solving \req{T_vars} for systems close to the kinetic equilibrium distribution \req{k_eq} but not necessarily in chemical equilibrium, a task for which the method in the previous section is not well suited in general. For a general kinetic equilibrium distribution, the temperature does not necessarily scale as $T\propto 1/a$ i.e. the temperature is not controlled solely by dilution.  For this reason, we will find it more useful to make the change of variables $z=p/T(t)$ rather than the scaling used in \req{a_vars}.  Here $T(t)$ is to be viewed as the time dependent effective temperature of the distribution $f$, a notion we will make precise later.  With this change of variables, the Boltzmann equation becomes
\begin{equation}\label{T_boltzmann}
\partial_t f-z\left(H+\frac{\dot T}{T}\right)\partial_z f=\frac{1}{E}C[f].
\end{equation}

 To model a distribution close to kinetic equilibrium at temperature $T$ and fugacity $\Upsilon$, we assume
\begin{equation}\label{kinetic_approx}
f(t,z)= f_\Upsilon (t,z)\psi(t,z)
\end{equation}
where the kinetic equilibrium distribution $f_\Upsilon $ depends on $t$ because we are assuming $\Upsilon$ is time dependent (with dynamics to be specified later). 

We will solve \req{T_boltzmann} by expanding $\psi$ in the basis of orthogonal polynomials, $\hat \psi_i$, generated by the parametrized weight function
\begin{equation}\label{weight}
w(z)\equiv w_\Upsilon(z)\equiv z^2f_\Upsilon (z)=\frac{z^2}{\Upsilon^{-1} e^z+1}
\end{equation}
on the interval $[0,\infty)$. See \ref{orthopoly_app} for details on the construction of these polynomials and their dependence on the parameter $\Upsilon$. This choice of weight is physically motivated by the fact that we are interested in solutions that describe massless particles not too far from kinetic equilibrium, but far from chemical equilibrium. We call this the chemical non-equilibrium method. 

We emphasize that we have made three important changes as compared to  the chemical equilibrium method:
\begin{enumerate}
\item  We allow a general time dependence of the effective temperature parameter $T$, i.e. we do not assume dilution temperature scaling $T=1/a$.
\item We have replaced the chemical  equilibrium distribution in the weight \req{free_stream_weight}  with a chemical non-equilibrium distribution  $f_\Upsilon $, i.e. we introduced $\Upsilon$.
\item We have introduced an additional factor of $z^2$ to the functional form of the weight as proposed in a different context in Refs.\cite{Wilkening,Wilkening2}. 
\end{enumerate} 
We note that the authors of \cite{Esposito2000} did consider the case of fixed chemical potential imposed as an initial condition. This is not the same as an emergent chemical non-equilibrium, i.e. time dependent $\Upsilon$, that we study here, nor do they consider a $z^2$ factor in the weight. We borrowed the idea for the $z^2$ prefactor from   Ref.\cite{Wilkening2}, where it was found that including a $z^2$ factor along with the non-relativistic chemical equilibrium distribution in the weight improved the accuracy of their method. Fortuitously,  this will also allow us to capture the particle number and energy flow with fewer terms than required by the chemical equilibrium method. A suitably modified weight and method allows us to maintain these advantages when a particle mass scale becomes relevant. We return to this problem in a subsequent work.


\subsection{ Comparison of Bases}\label{basis_comparison}

\begin{figure}[H]
 \begin{minipage}[b]{0.5\linewidth}
\centerline{\includegraphics[height=6.cm]{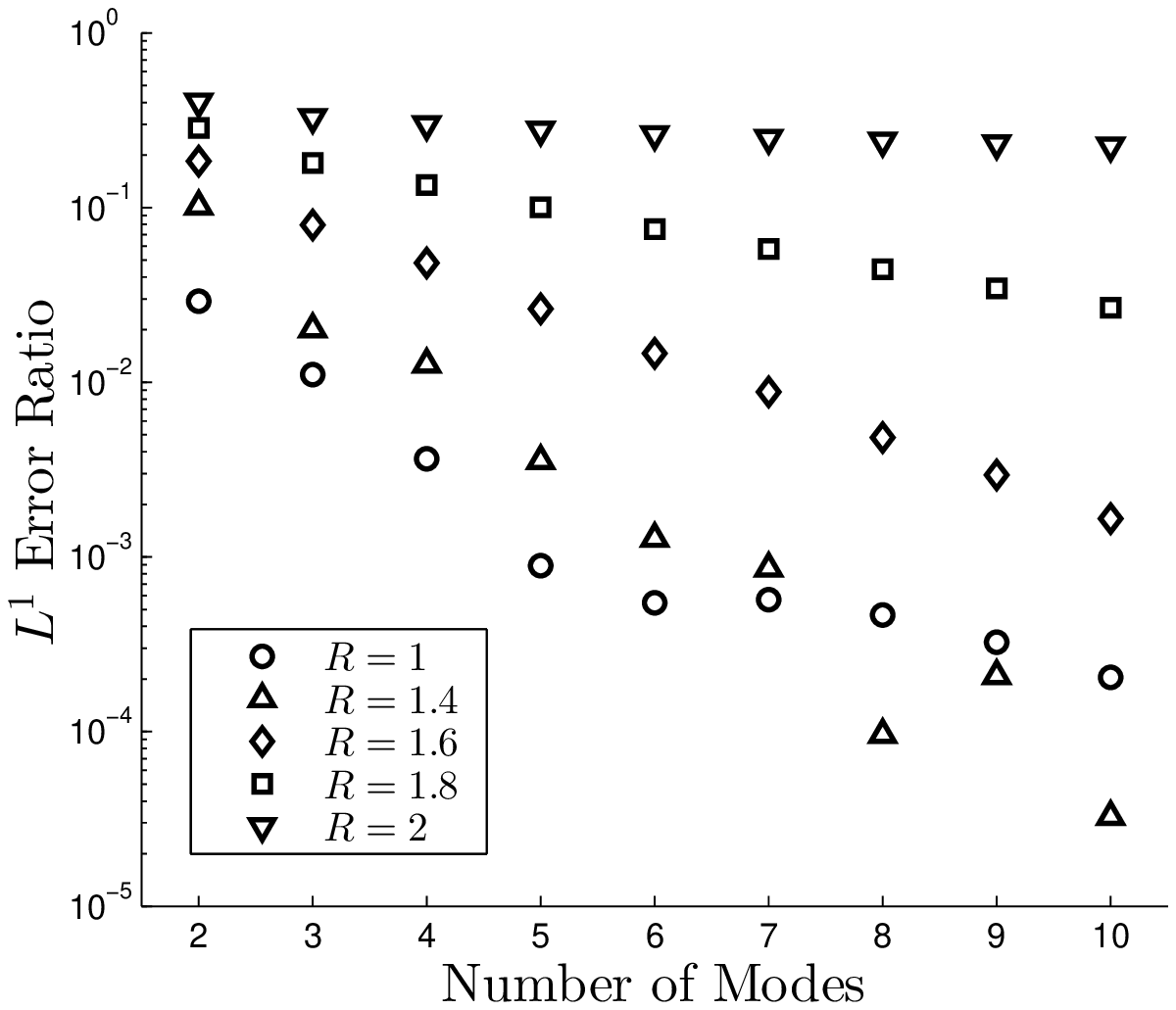}}
\caption{Errors in expansion of \req{zeroth_approx} as a function of number of modes, $\Upsilon=0.5$.}\label{fig:free_stream_f0_approx_Ups_5}
 \end{minipage}
 \hspace{0.5cm}
 \begin{minipage}[b]{0.5\linewidth}
\centerline{\includegraphics[height=6cm]{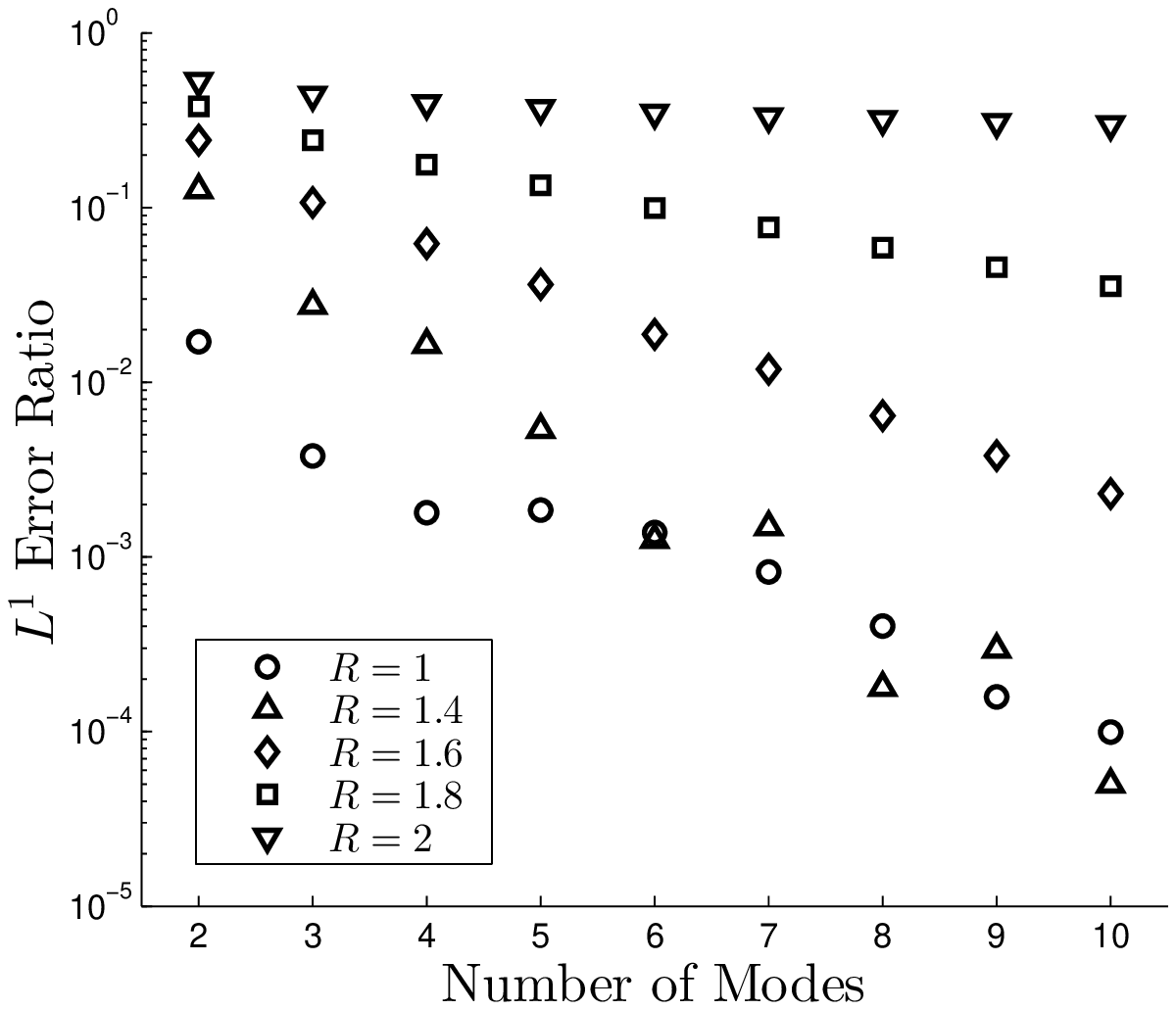}}
\caption{Errors in  expansion of \req{zeroth_approx} as a function of number of modes, $\Upsilon=1.5$.}\label{fig:free_stream_f0_approx_Ups_1_5}
 \end{minipage}
 \end{figure}

Before deriving the dynamical equations for the method outlined in section \ref{kinetic_eq_approach}, we illustrate the error inherent in approximating the chemical non-equilibrium distribution \req{k_eq}  with a  chemical equilibrium distribution \req{ch_eq} whose temperature is $T=1/a$.   Given a chemical non-equilibrium distribution 
\begin{equation}\label{zeroth_approx}
f_\Upsilon (y)=\frac{1}{\Upsilon^{-1}e^{y/(aT)}+1},
\end{equation}
 we can attempt to write it as a perturbation of the chemical equilibrium distribution,  
\begin{equation}\label{chi_def}
f_\Upsilon=f_{ch}\chi
\end{equation} as we would need to when using the method of section \ref{free_stream_approach}.  We expand $\chi=\sum_i d^i\hat\chi_i$ in the orthonormal basis generated by $f_{ch}$ and, using $N$ terms, form the $N$-mode approximation $f_\Upsilon^N$ to $f_\Upsilon$.  The $d^i$ are obtained by taking the $L^2(f_{ch}dy)$ inner product of $\chi$ with the basis function $\hat\chi_i$,
\begin{equation}
d^i=\int\hat\chi_i \chi f_{ch}dy=\int\hat\chi_i  f_\Upsilon dy.
\end{equation}
 Figures \ref{fig:free_stream_f0_approx_Ups_5} and \ref{fig:free_stream_f0_approx_Ups_1_5} show the normalized $L^1(dx)$ errors between $f_\Upsilon^N$ and $f_\Upsilon$, computed via
\begin{equation}
\text{error}_N=\frac{\int_0^\infty |f_\Upsilon -f_\Upsilon ^N|dy}{\int_0^\infty |f_\Upsilon |dy}.
\end{equation}

We note the appearance of the reheating ratio
\begin{equation}\label{reheat}
 R\equiv aT  
\end{equation}
in the denominator of \req{zeroth_approx}, which comes from changing variables from $z=p/T$ in \req{weight} to $y=ap$ in order to compare with \req{free_stream_weight}.  Physically, $R$ is the ratio of the physical temperature $T$ to the dilution controlled temperature scaling  of $1/a$.   In physical situations, including cosmology, $R$ can vary from unity when dimensioned energy scales influence dynamical equations for $a$. From the error plots we see that for $R$ sufficiently close to $1$, the approximation performs well with a small number of terms, even with $\Upsilon\neq 1$.  

\begin{figure}[H]
 \begin{minipage}[b]{0.5\linewidth}
\centerline{\includegraphics[height=6.cm]{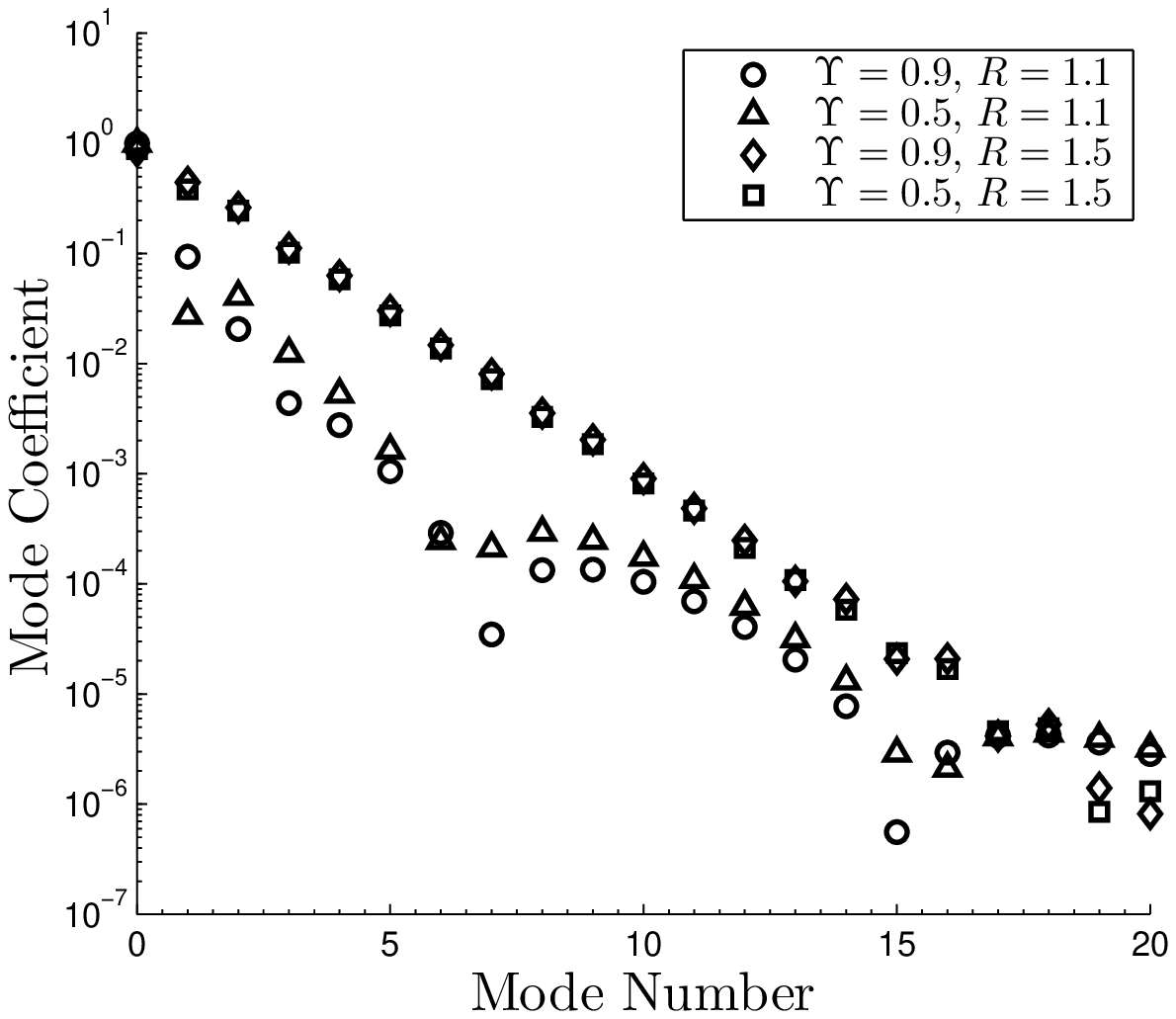}}
\caption{Mode coefficients for the expansion of $\hat f_{\Upsilon}$ in terms of the chemical equilbrium basis.}\label{fig:free_stream_mode_coeff}
 \end{minipage}
 \hspace{0.5cm}
 \begin{minipage}[b]{0.5\linewidth}
\centerline{\includegraphics[height=6cm]{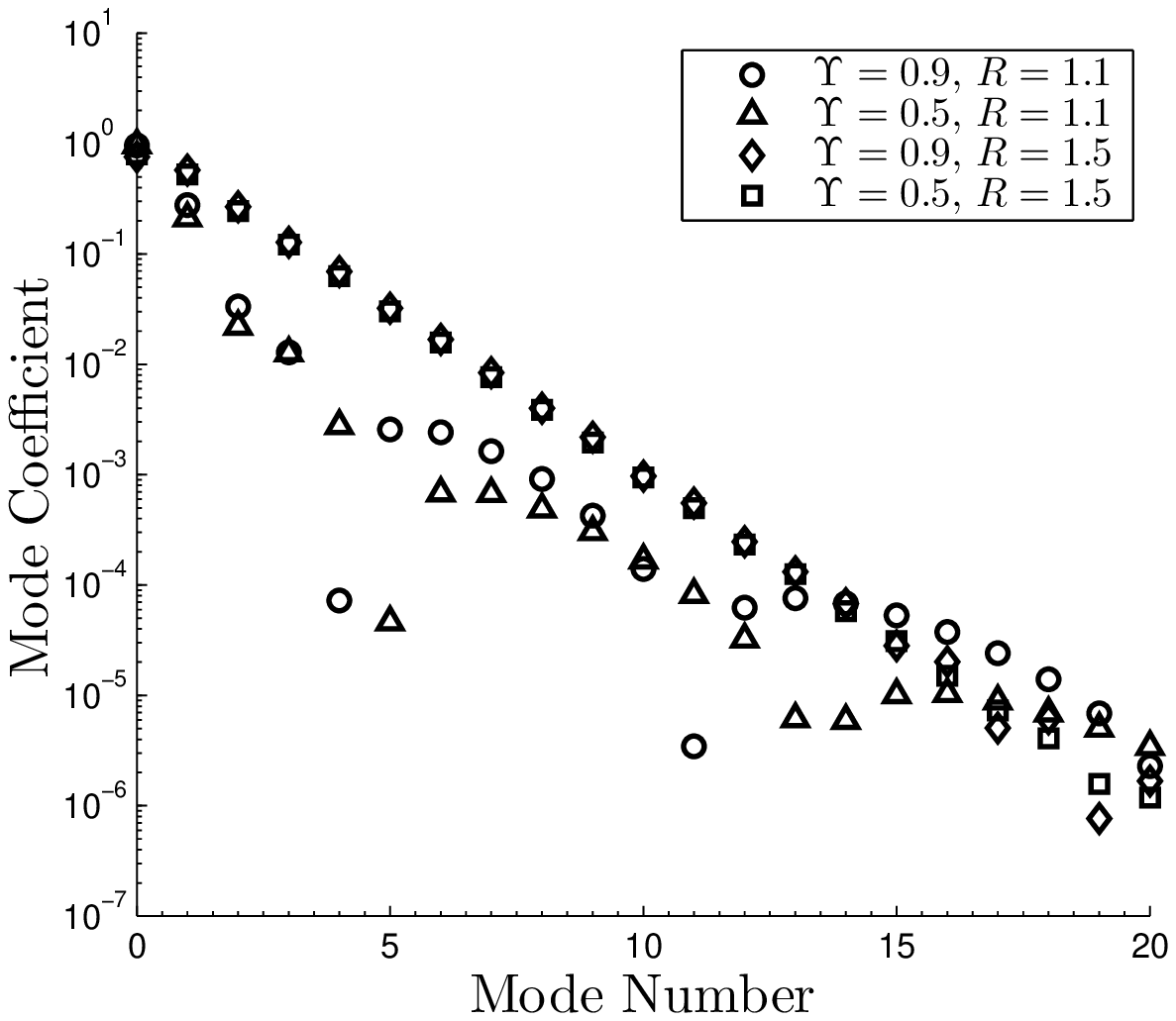}}
\caption{Mode coefficients for the expansion of $\hat f_{\Upsilon}$ in terms of the basis of Laguerre polynomials.}\label{fig:laguerre_mode_coeff}
 \end{minipage}
 \end{figure}

 For non-degenerate systems (i.e. $\Upsilon$ small) the chemical non-equilibrium weight is close to the weight that generates the Laguerre polynomials, $e^{-y}$. In figures \ref{fig:free_stream_mode_coeff} and \ref{fig:laguerre_mode_coeff} we show the mode coefficients for the expansion of $\hat f_\Upsilon$ in term of the  the basis generated by $f_{ch}$ and Laguerre polynomials respectively. The hat denotes that we divide by the corresponding $L^2$ norm of $f_\Upsilon$ in order to properly compare the two spaces.

 For $\Upsilon\approx 1$ and $R\approx 1$ we have $f_{ch}\approx f_\Upsilon$ and so one expects the chemical equilibrium basis to provide a better approximation than Laguerre polynomials.  The mode coefficients are significantly higher in figure \ref{fig:laguerre_mode_coeff}  as compared to  figure \ref{fig:free_stream_mode_coeff} for the case $\Upsilon=0.9$, $R=1.1$. This indicates a larger truncation error and hence a worse approximation using Laguerre polynomials, as expected.  We note that the regime in which the chemical equilibrium method was employed in \cite{Esposito2000,Mangano2002} was very close to $\Upsilon=1$, $R=1$, making their use of the chemical equilibrium method an appropriate choice.

 In contrast, when the system is non-degenerate or when there is significant reheating we see that there is little to no advantage in using the chemical equilibrium basis over the basis of Laguerre polynomials, as seen by the similar truncation errors between the two plots when $R=1.5$ and, to a slighly lesser degree, $\Upsilon=0.5$.  It is in these situations where the chemical non-equilibrium method we present becomes most advantageous, being superior to both the Laguerre and chemical equilibrium methods.

\begin{figure}[H]
 \begin{minipage}[b]{0.5\linewidth}
\centerline{\includegraphics[height=6.cm]{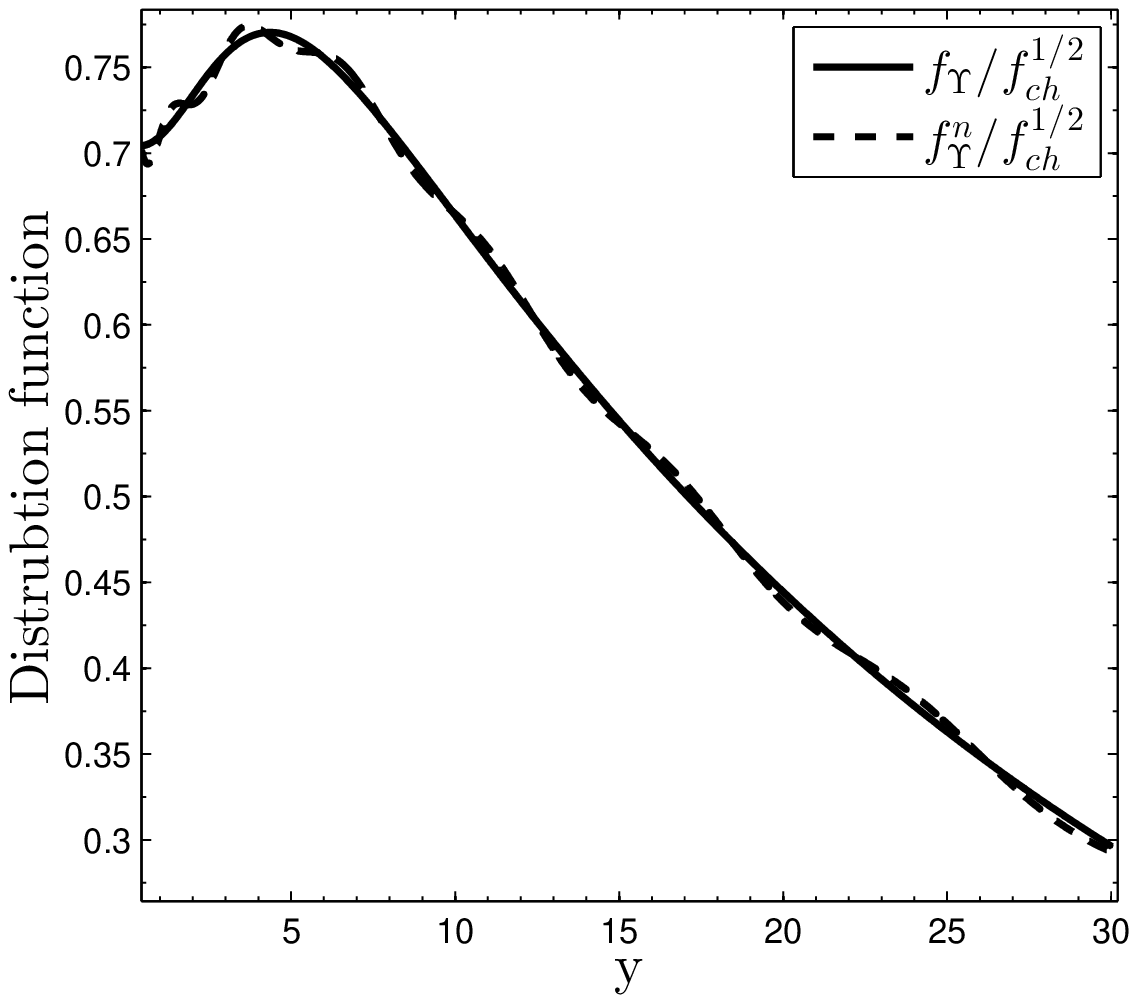}}
\caption{Approximation to \req{zeroth_approx} for $\Upsilon=1$ and $R=1.85$ using the first $20$ basis elements generated by \req{free_stream_weight}.}\label{fig:free_stream_f0_approx_Ups_1_T_r_1_85}
 \end{minipage}
 \hspace{0.5cm}
 \begin{minipage}[b]{0.5\linewidth}
\centerline{\includegraphics[height=6.cm]{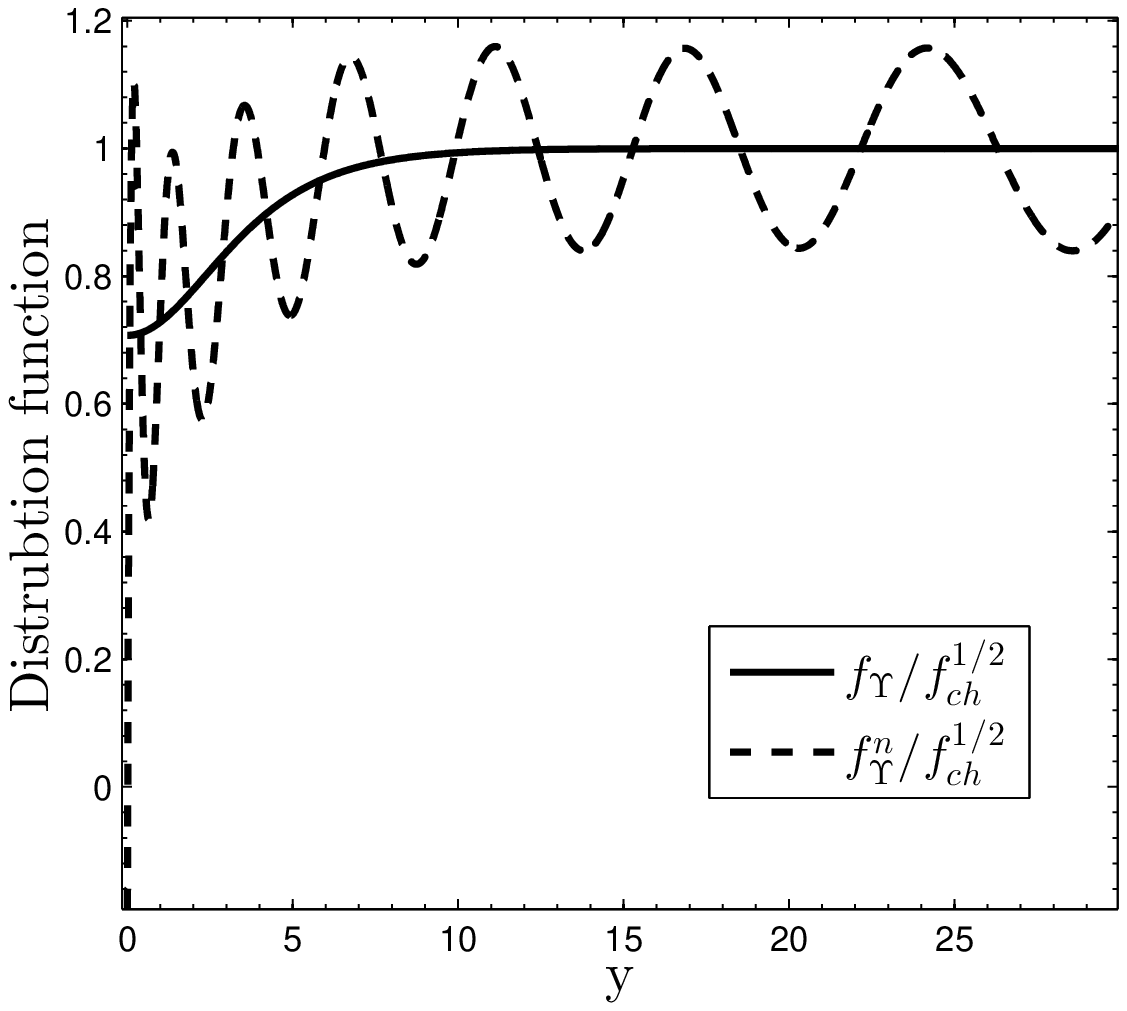}}
\caption{Approximation to \req{zeroth_approx} for $\Upsilon=1$ and $R=2$ using the first $20$ basis elements generated by \req{free_stream_weight}.}\label{fig:free_stream_f0_approx_Ups_1_T_r_2}
 \end{minipage}
 \end{figure}
 In the case of large reheating, we find that when $R$ approaches and surpasses $2$, large spurious oscillations begin to appear in the expansion and they persist even when a large number of terms are used, as seen in figures  \ref{fig:free_stream_f0_approx_Ups_1_T_r_1_85} and \ref{fig:free_stream_f0_approx_Ups_1_T_r_2}, where we compare $f_\Upsilon/f_{ch}^{1/2}$ with $f_{\Upsilon}^N/f_{ch}^{1/2}$ for $\Upsilon=1$ and $N=20$.  The poor behavior arises from the fact that for $R\geq 2$,  the perturbation $\chi$ from \req{chi_def} is no longer contained in the Hilbert space $L^2(f_{ch}dy)$.  This can be seen as follows.  The squared norm of $\chi$ is
\begin{equation}
\int\chi^2f_{ch}dy=\int\frac{f_\Upsilon^2}{f_{ch}}dy=\int\frac{e^y+1}{(\Upsilon^{-1} e^{y/R}+1)^2}dy.
\end{equation}
For $y$ large, the integrand can be approximated by 
\begin{equation}
\frac{e^y+1}{(\Upsilon^{-1} e^{y/R}+1)^2}\approx \Upsilon^2e^{y(1-2/R)}
\end{equation}
which is not integrable for $R\geq 2$.

This demonstrates that the chemical equilibrium method with dilution temperature scaling will  perform  poorly in situations that experience a large degree of reheating (and so will Laguerre polynomials for the same reason). The dynamical effective temperature that we introduce in the chemical non-equilibrium method will allow us to handle such cases without issue.

For $R\approx 1$, the benefit of including fugacity is not as striking, as the chemical equilibrium basis is able to approximate \req{zeroth_approx} reasonably well.  However, for more stringent error tolerances including $\Upsilon$ can reduce the number of required modes in cases where the degree of chemical non-equilibrium is large.

\subsection{Dynamics}\label{dynamics_sec}
In this section we derive the dynamical equations for the  method outlined in section \ref{kinetic_eq_approach}.  In particular, we identify physically motivated dynamics for the effective temperature and fugacity.  Using \req{T_boltzmann} and the definition of $\psi$ from \req{kinetic_approx} we have
\begin{align}\label{near_equilib_eq}
\partial_t \psi+\frac{1}{f_\Upsilon }\frac{\partial f_\Upsilon }{\partial\Upsilon}\dot\Upsilon\psi-\frac{z}{f_\Upsilon }\left(H+\frac{\dot{T}}{T}\right)\left(\psi\partial_zf_\Upsilon +f_\Upsilon \partial_z \psi\right)=\frac{1}{f_\Upsilon E}C[f_\Upsilon \psi].
\end{align}
 Denote the monic orthogonal polynomial basis generated by the weight \req{weight} by $\psi_n$, $n=0,1,...$ where $\psi_n$ is degree $n$ and call the normalized versions  $\hat{\psi}_n$. Recall that $\hat\psi_n$ depend on $t$ due to the $\Upsilon$ dependence of the weight function used in the construction. Consider the space of polynomial of degree less than or equal to $N$, spanned by $\hat\psi_n$, $n=0,...,N$.   For $\psi$ in this subspace, we expand $\psi=\sum_{j=0}^Nb^j\hat\psi_j$ and use \req{near_equilib_eq}  to obtain
\begin{align}\label{T_vars}
\sum_i \dot{b}^i\hat\psi_i=&\sum_ib^i\frac{z}{f_\Upsilon }\left(H+\frac{\dot{T}}{T}\right)\left(\partial_z(f_\Upsilon )\hat\psi_i+f_\Upsilon \partial_z\hat\psi_i\right)\\
&-\sum_ib^i\left(\dot{\hat{\psi}}_i+\frac{1}{f_\Upsilon }\frac{\partial f_\Upsilon }{\partial\Upsilon}\dot\Upsilon\hat\psi_i\right)+\frac{1}{f_\Upsilon E}C[f].
\notag
\end{align}
From this we see  that the equations obtained from the Boltzmann equation by projecting onto the finite dimensional subspace are
\small
\begin{align}
\dot b^k=& \sum_i b^i\left(H+\frac{\dot{T}}{T}\right)\left(\langle\frac{z}{f_\Upsilon }\hat\psi_i\partial_zf_\Upsilon ,\hat\psi_k\rangle+\langle z\partial_z \hat\psi_i,\hat\psi_k\rangle\right) \\
&-\sum_i b^i\dot{\Upsilon}\left(\langle\frac{1}{f_\Upsilon }\frac{\partial f_\Upsilon }{\partial\Upsilon}\hat\psi_i,\hat\psi_k\rangle+\langle\frac{\partial\hat{\psi}_i}{\partial \Upsilon},\hat\psi_k\rangle\right)+\langle\frac{1}{f_\Upsilon E}C[f],\hat\psi_k\rangle\notag
\end{align}
\normalsize
where $\langle\cdot,\cdot\rangle$ denotes the inner product defined by the weight function \req{weight}
\begin{equation}
\langle h_1,h_2\rangle=\int_0^\infty h_1(z)h_2(z)w_\Upsilon(z)dz.
\end{equation}
  The term in brackets comprises the linear part of the system, while the collision term contains polynomial nonlinearities when multiple coupled distribution are being modeled using a $2$-$2$ collision operator \req{coll}.  

To isolate the linear part, we define matrices
\begin{align}\label{A_B_matrices}
A^k_i(\Upsilon)\equiv&\langle\frac{z}{f_\Upsilon }\hat\psi_i\partial_zf_\Upsilon ,\hat\psi_k\rangle+\langle z\partial_z \hat\psi_i,\hat\psi_k\rangle,\\
B^k_i(\Upsilon)\equiv& C_i^k(\Upsilon)+D_i^k(\Upsilon),\hspace{2mm} C_i^k\equiv\Upsilon\langle\frac{1}{f_\Upsilon }\frac{\partial f_\Upsilon }{\partial\Upsilon}\hat\psi_i,\hat\psi_k\rangle,\hspace{2mm} D_i^k\equiv\Upsilon\langle\frac{\partial\hat{\psi}_i}{\partial \Upsilon},\hat\psi_k\rangle. 
\end{align}
 With these definitions, the equations for the $b^k$ become
\begin{align}\label{b_eq}
\dot b^k=& \left(H+\frac{\dot{T}}{T}\right)\sum_i A_i^k(\Upsilon)b^i-\frac{\dot{\Upsilon}}{\Upsilon}\sum_i B_i^k(\Upsilon)b^i+\langle\frac{1}{f_\Upsilon E}C[f],\hat\psi_k\rangle.
\end{align}
 See \ref{ortho-general} for details on how to recursively construct the $\partial_z\hat\psi_i$. We show how to compute the inner products $\langle\hat\psi_k,\partial_{\Upsilon}\hat\psi_k\rangle$ in  \ref{ortho-polynom-fam}. In \ref{lower_triang} we prove that that both $A$ and $B$ are lower triangular and show that the only inner products involving the $\partial_\Upsilon\hat{\psi}_i$ that are required in order to compute $A$ and $B$ are those the above mentioned diagonal elements, $\langle\hat\psi_k,\partial_{\Upsilon}\hat\psi_k\rangle$.

We fix the dynamics of $T$ and $\Upsilon$ by imposing the conditions
\begin{equation}\label{b_ics}
b^0(t)\hat\psi_0(t)=1,\hspace{2mm}b^1(t)=0.
\end{equation}
In other words,
\begin{equation}
f(t,z)=f_\Upsilon (t,z)\left(1+\phi(t,z)\right),\hspace{2mm} \phi=\sum_{i=2}^N b^i\hat\psi_i.
\end{equation}
This reduces the number of degrees of freedom in \req{b_eq} from $N+3$ to $N+1$.  In other words, after enforcing \req{b_ics}, \req{b_eq} constitutes $N+1$ equations for the remaining $N+1$ unknowns, $b^2,...,b^N$, $\Upsilon$, and $T$.  We will call $T$ and $\Upsilon$ the first two ``modes", as their dynamics arise from imposing the conditions \req{b_ics} on the zeroth and first order coefficients in the expansion. We will solve for their dynamics explicitly below.

To see the physical motivation for the choices \req{b_ics}, consider the particle number density and energy density.  Using orthonormality of the $\hat\psi_i$ and \req{b_ics} we have
\begin{align}
n=&\frac{g_pT^3}{2\pi^2}\sum_ib^i\int_0^\infty f_\Upsilon  \hat\psi_i z^2 dz=\frac{g_pT^3}{2\pi^2}\sum_ib^i\langle \hat\psi_i ,1\rangle\\
=&\frac{g_pT^3}{2\pi^2}b^0\langle \hat\psi_0 ,1\rangle=\frac{g_pT^3}{2\pi^2}\langle 1 ,1\rangle,\\
\rho=&\frac{g_pT^4}{2\pi^2}\sum_ib^i\int_0^\infty f_\Upsilon  \hat\psi_i z^3 dz=\frac{g_pT^4}{2\pi^2}\sum_ib^i\langle\hat\psi_i, z\rangle\\
=&\frac{g_pT^4}{2\pi^2}\left(b^0\langle\hat\psi_0, z\rangle+b^1\langle\hat\psi_1, z\rangle\right)=
\frac{g_pT^4}{2\pi^2}\langle 1,z\rangle.
\end{align}
 Using these together with the definition of the weight function \req{weight} we find
\begin{align}\label{th_eq_moments}
n=&\frac{g_pT^3}{2\pi^2}\int_0^\infty f_\Upsilon  z^2dz,\\
\label{th_eq_moments2}
\rho=&\frac{g_pT^4}{2\pi^2}\int_0^\infty f_\Upsilon  z^3dz.
\end{align}

Equations (\ref{th_eq_moments}) and (\ref{th_eq_moments2}) show that the first two modes, $T$ and $\Upsilon$, with time evolution fixed by \req{b_ics} combine with the chemical non-equilibrium distribution $f_\Upsilon $ to capture the particle number density and energy density of the system exactly.  This fact is very significant, as it implies that within the chemical non-equilibrium approach as long as the back-reaction from the non-thermal distortions is small (meaning that the evolution of $T(t)$ and $\Upsilon(t)$ is not changed significantly when more modes are included), {\em all the effects relevant to the computation of  particle and energy flow are modeled by the time evolution of $T$ and $\Upsilon$ alone} and no further modes are necessary.  This gives a clear separation between the averaged physical quantities, characterizing the time evolution of $f_\Upsilon $, and the momentum dependent non-thermal distortions as contained in 
\begin{equation}
\phi=\sum_{i=2}^N b^i\hat\psi_i.
\end{equation}

One should contrast this chemical non-equilibrium behavior  with the chemical equilibrium situation, where a minimum of four modes is required to describe the number and energy densities, as shown in \req{free_stream_moments}.   Moreover we will show that convergence to the desired precision is faster in the chemical non-equilibrium approach as compared to the chemical equilibrium. Due to the high cost of numerically integrating realistic collision integrals of the form \req{coll}, this fact can be very significant in applications. We remark that the relations \req{th_eq_moments} are the physical motivation for including the $z^2$ factor in the weight function. All three modifications we have made in constructing our new method, the introduction of an effective temperature i.e. $R\ne 1$, the generalization to chemical non-equilibrium $f_\Upsilon $, and the introduction of $z^2$ to the weight, \req{reheat}, were needed to obtain the properties \req{th_eq_moments}, but it is the introduction of $z^2$ that reduces the number of required modes and hence reduces the computational cost. 

With $b^0$ and $b^1$ fixed as in \req{b_ics} we can solve the equations for $\dot b^0$ and $\dot b^1$ from \req{b_eq} for $\dot\Upsilon$ and $\dot T$ to obtain
\small
\begin{align}\label{Ups_T_eqs}
\dot{\Upsilon}/{\Upsilon}=&\frac{(Ab)^1\langle\frac{1}{f_\Upsilon E}C[f],\hat\psi_0\rangle-(Ab)^0\langle\frac{1}{f_\Upsilon E}C[f],\hat\psi_1\rangle }{[\Upsilon\partial_\Upsilon \langle1,1\rangle/(2||\psi_0||)+(Bb)^0](Ab)^1-(Ab)^0(Bb)^1},\\[0.5cm]
\dot{T}/T
=&\frac{(Bb)^1\langle\frac{1}{f_\Upsilon E}C[f],\hat\psi_0\rangle-\langle\frac{1}{f_\Upsilon E}C[f],\hat\psi_1\rangle[\Upsilon\partial_\Upsilon \langle1,1\rangle/(2||\psi_0||)+(Bb)^0]}{[\Upsilon\partial_\Upsilon \langle1,1\rangle/(2||\psi_0||)+(Bb)^0](Ab)^1-(Ab)^0(Bb)^1}-H\notag\\[0.3cm]
=&\frac{1}{(Ab)^1}\left((Bb)^1\dot{\Upsilon}/\Upsilon-\langle\frac{1}{f_\Upsilon E}C[f],\hat\psi_1\rangle\right)-H.\label{T_eq}
\end{align}
\normalsize
Here $(Ab)^n=\sum_{j=0}^NA^n_jb^j$ and similarly for $B$ and $||\cdot||$ is the norm induced by $\langle\cdot,\cdot\rangle$. In deriving this, we used
\begin{equation}
\dot{b}^0=\frac{1}{2||\psi_0||}\dot\Upsilon\partial_\Upsilon \langle1,1\rangle, \hspace{4mm} \partial_\Upsilon \langle1,1\rangle=\int_0^\infty \frac{z^2}{(e^{z/2}+ \Upsilon e^{-z/2})^2}dz
\end{equation}
which comes from differentiating \req{b_ics}. 
 
 It is easy to check that when the collision operator vanishes, then the above system is solved by 
\begin{equation}\label{free_stream_sol}
\Upsilon=\text{constant},\hspace{4mm} \frac{\dot T}{T}=-H,\hspace{2mm}  b^n=\text{constant},\hspace{1mm} n>2
\end{equation}
i.e. the fugacity and non-thermal distortions are `frozen' into the distribution and the temperature satisfies dilution scaling $T\propto 1/a$.

When the collision term becomes small, \req{free_stream_sol} motivates another change of variables. Letting $T=(1+\epsilon)/a$  gives the equation
\begin{equation}
\dot\epsilon=\frac{1+\epsilon}{(Ab)^1}\left((Bb)^1\dot{\Upsilon}/\Upsilon-\langle\frac{1}{f_\Upsilon E}C[f],\hat\psi_1\rangle\right).
\end{equation}
Solving this in place of \req{T_eq} when the collision terms are small avoids having to numerically track the free-streaming evolution.  In particular this will ensure conservation of comoving particle number, which equals a function of $\Upsilon$ multiplied by $(aT)^3$, to much greater precision in this regime as well as resolve the freeze-out temperatures more accurately.

\subsubsection{Projected Dynamics are Well-defined}
The following calculation shows that, for a distribution initially in kinetic equilibrium, the determinant factor in the denominator of \req{Ups_T_eqs} is nonzero and hence the dynamics for $T$ and $\Upsilon$, as well as the remainder of the projected system, are well-defined, at least for sufficiently small times. 

 Kinetic equilibrium implies the initial conditions $b^0=||\psi_0||$, $b^i=0$, $i>0$. Defining
\begin{equation}
K\equiv (\Upsilon\partial_\Upsilon \langle 1,1\rangle/(2||\psi_0||)+(Bb)^0)(Ab)^1-(Ab)^0(Bb)^1\
\end{equation}
and using the definition \req{A_B_matrices} for $B$ we have
\begin{align}
K=&(C^0_0A^1_0-A^0_0C^1_0)(b^0)^2+\left[(D^0_0A^1_0-A^0_0D^1_0)(b^0)^2+\Upsilon\partial_\Upsilon \langle 1,1\rangle/(2||\psi_0||)A^1_0b^0\right]\notag\\[0.3cm]
\equiv & K_1+K_2.\notag
\end{align}
 Using the definitions of $C$ and $D$ from  \req{A_B_matrices}  we find
\begin{align}
K_1=&\langle \frac{1}{1+\Upsilon e^{-z}},1\rangle\langle \frac{-z}{1+\Upsilon e^{-z}}\hat\psi_1,\hat\psi_0\rangle-\langle\frac{-z}{1+\Upsilon e^{-z}},\hat\psi_0\rangle\langle\frac{1}{1+ \Upsilon e^{-z}}\hat\psi_1,1\rangle.
\end{align}
Inserting the formula for $\hat\psi_1$ from \req{poly_recursion} we find
\begin{align}
K_1=&-\frac{1}{||\psi_1||\,||\psi_0||}\left[\langle\frac{1}{1+ \Upsilon e^{-z}},\hat\psi_0\rangle\langle\frac{z^2}{1+\Upsilon e^{-z}},\hat\psi_0\rangle-\langle\frac{z}{1+\Upsilon e^{-z}},\hat\psi_0\rangle^2\right].
\end{align}
The Cauchy-Schwarz inequality  applied to the inner product with weight function
\begin{equation}
\tilde{w}=\frac{w}{1+\Upsilon e^{-z}}\hat\psi_0
\end{equation}
together with linear independence of $1$ and $z$ implies that the term in brackets is positive and so $K_1<0$ at $t=0$.  For the second term, noting that $D^1_0=0$ by orthogonality and using \req{norm_deriv_eq}, we have

\begin{align}
K_2=&[\langle\partial_\Upsilon\hat\psi_0,\hat\psi_0\rangle||\psi_0||+\partial_\Upsilon \langle 1,1\rangle/(2||\psi_0||)]\Upsilon A_0^1||\psi_0||\\
=&0.\notag
\end{align}
This proves that $K$ is nonzero at $t=0$.\\

\section{Validation}\label{validation}
We will validate our numerical method on an exactly solvable model problem
\begin{equation}\label{toy_eq}
\partial_t f-pH \partial_p f=M\left(\frac{1}{\Upsilon^{-1}e^{p/T_{eq}}+1}-f(p,t)\right), \hspace{2mm} f(p,0)=\frac{1}{e^{p/T_{eq}(0)}+1}
\end{equation}
where $M$ is a constant with units of energy and we choose units in which it is equal to $1$. This model describes a distribution that is attracted to a given equilibrium distribution at a prescribed time dependent temperature $T_{eq}(t)$ and fugacity $\Upsilon$. This type of an idealized scattering operator, without fugacity, was first introduced in \cite{Anderson_Witting}.  For an application of this method involving more complex collision terms of the form in \req{coll} see \cite{Birrell_nu_param}.

 By changing coordinates $y=a(t)p$ we can rewrite \req{toy_eq} as
\begin{equation}\label{free_stream_toy}
\partial_tf(y,t)=\frac{1}{\Upsilon^{-1}\exp[y/(a(t)T_{eq}(t))]+1}-f(y,t).
\end{equation}
 which has as solution
\begin{equation}
f(y,t)=\int_0^t\frac{e^{s-t}}{\Upsilon^{-1}\exp[y/(a(s)T_{eq}(s))]+1}ds+\frac{e^{-t}}{\exp[y/(a(0)T_{eq}(0))]+1}.
\end{equation}
We now transform to $z=p/T(t)$ where the temperature $T$ of the distribution $f$ is defined as in section \ref{dynamics_sec}. Recalling that $H(t)=\dot a/a$, we have the exact solution to 
\begin{equation}\label{k_eq_toy}
\partial_tf(t,z)-z\left(H+\frac{\dot{T}}{T}\right)\partial_zf(t,z)=\frac{1}{\Upsilon^{-1}e^{zT/T_{eq}}+1}-f(z,t)
\end{equation}
given by
\begin{align}\label{exact_sol}
f(z,t)=&\int_0^t\frac{e^{s-t}}{\Upsilon^{-1}\exp[a(t)T(t)z/(a(s)T_{eq}(s))]+1}ds\\
&+\frac{e^{-t}}{\exp[a(t)T(t)z/(a(0)T_{eq}(0))]+1}.\notag
\end{align}
We use this to test the chemical equilibrium and chemical non-equilibrium methods under two different conditions. 

\subsection{Reheating Test}
First we test the two methods we have outlined in a scenario that exhibits reheating.  Motivated by applications to cosmology, we choose a scale factor evolving as in the radiation dominated era, a fugacity $\Upsilon=1$, and choose an equilibrium temperature that exhibits reheating like behavior with $aT_{eq}$ increasing for a period of time,
\begin{align}\label{a_T_def}
a(t)=\left(\frac{t+b}{b}\right)^{1/2}\!\!\!,\ \  \ \
T_{eq}(t)=\frac{1}{a(t)}\left(1+\frac{1-e^{-t}}{e^{-(t-b)}+1}(R-1)\right)
\end{align}
where $R$ is the desired reheating ratio. Note that $(aT_{eq})(0)=1$ and $(aT_{eq})(t)\rightarrow R$ as $t\rightarrow\infty$. Qualitatively, this is reminiscent of the dynamics of neutrino freeze-out, but the range of reheating ratio for which we will test our method is larger than found there.

We solved \req{free_stream_toy} and \req{k_eq_toy} numerically using the chemical equilibrium and chemical non-equilibrium methods respectively for $t\in[0,10]$ and $b=5$ and the cases $R=1.1$, $R=1.4$, as well as the more extreme ratio of $R=2$.  The bases of orthogonal polynomials were generated numerically using the recursion relations from \ref{orthopoly_app}.  For the applications we are considering, where the solution is a small perturbation of equilibrium, only a small number of terms are required and so the numerical challenges associated with generating a large number of such orthogonal polynomials are not an issue.

\subsubsection{Chemical Equilibrium Method}
We solved \req{free_stream_toy} using the chemical equilibrium method, with the orthonormal basis defined by the weight function \req{free_stream_weight} for $N=2,...,10$ modes (mode numbers $n=0,...,N-1$) and prescribed  single step relative and absolute error tolerances of $10^{-13}$ for the numerical integration, and with asymptotic reheating ratios of $R=1.1$, $R=1.4$, and $R=2$.

\begin{figure}[H]
 \begin{minipage}[b]{0.5\linewidth}
\centerline{\includegraphics[height=6.2cm]{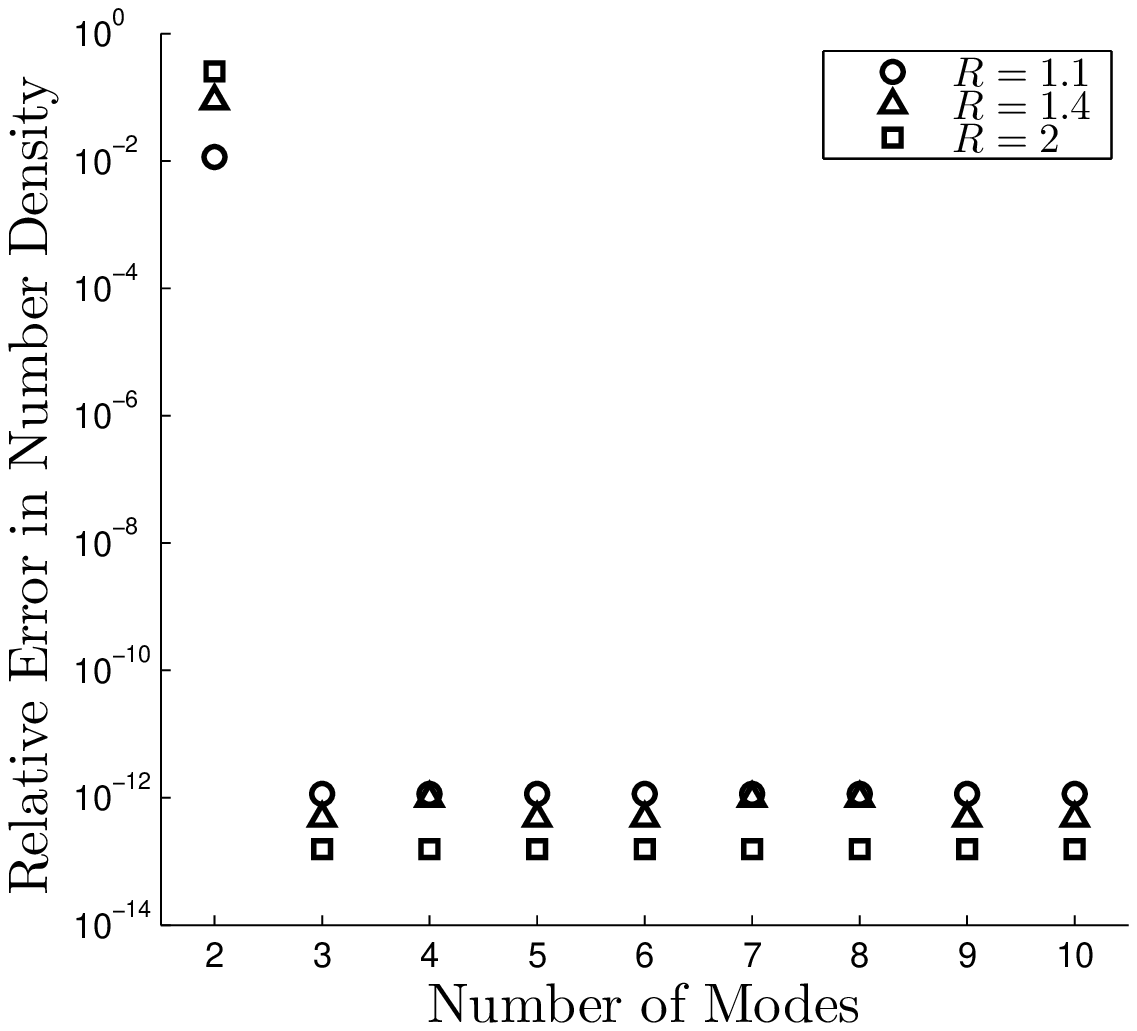}}
\caption{Maximum relative error in particle number density.}\label{fig:free_stream_num_err}
 \end{minipage}
 \hspace{0.5cm}
 \begin{minipage}[b]{0.5\linewidth}
\centerline{\includegraphics[height=6.2cm]{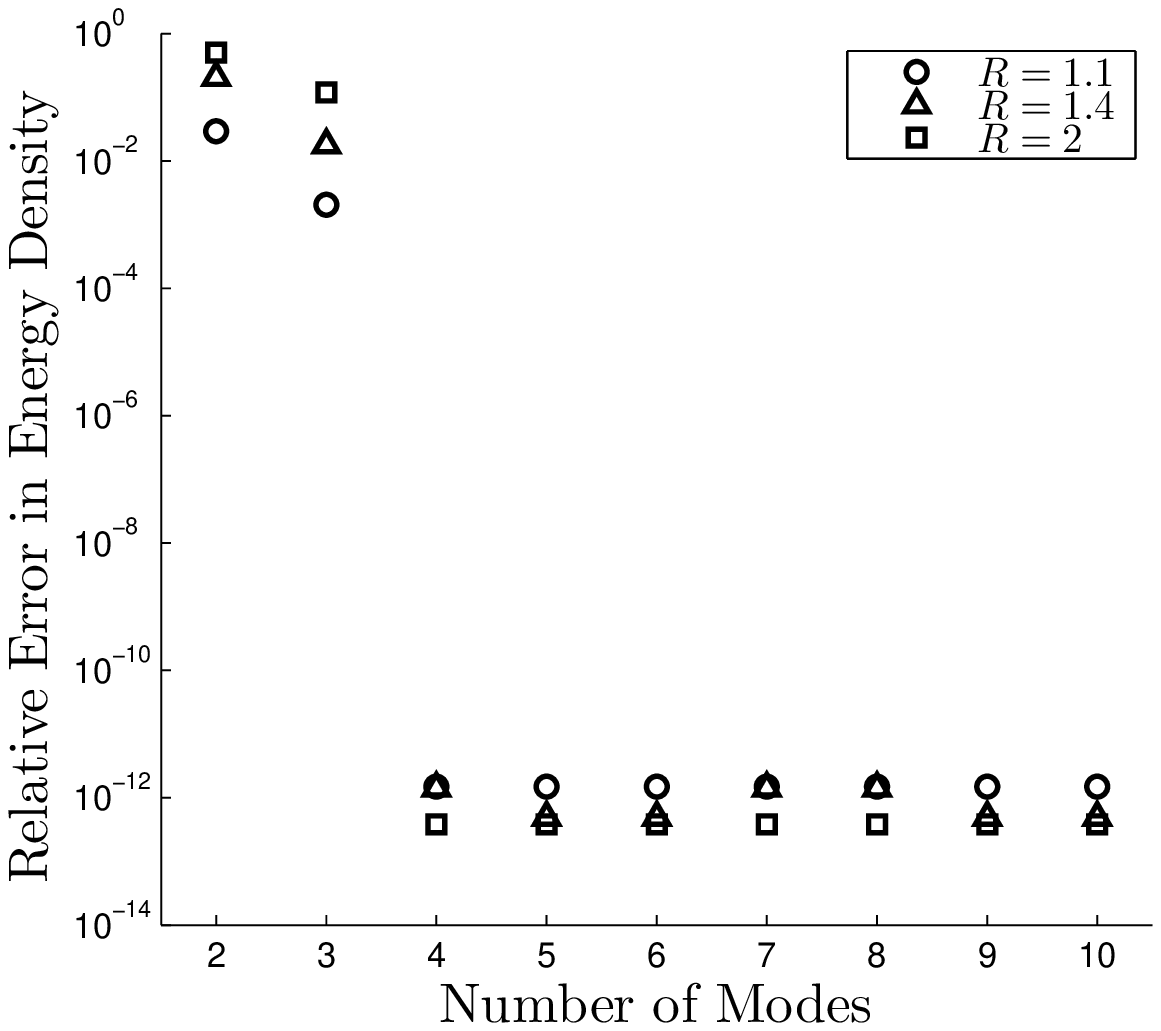}}
\caption{Maximum relative error in energy density.}\label{fig:free_stream_E_err}
 \end{minipage}
 \end{figure}

In figures   \ref{fig:free_stream_num_err} and  \ref{fig:free_stream_E_err} we show the maximum relative error in the number densities and energy densities respectively over the time interval $[0,10]$ for various numbers of computed modes.  The particle number density and energy density are accurate, up to the integration tolerance level, for $3$ or more and $4$ or more modes respectively. This is consistent with \req{free_stream_moments} which shows the number of modes required to capture each of these quantities. However, fewer modes than these minimum values lead to a large error in the corresponding moment of the distribution function.

\begin{figure}[H]
\begin{minipage}[t]{0.5\linewidth}
\centerline{\includegraphics[height=6cm]{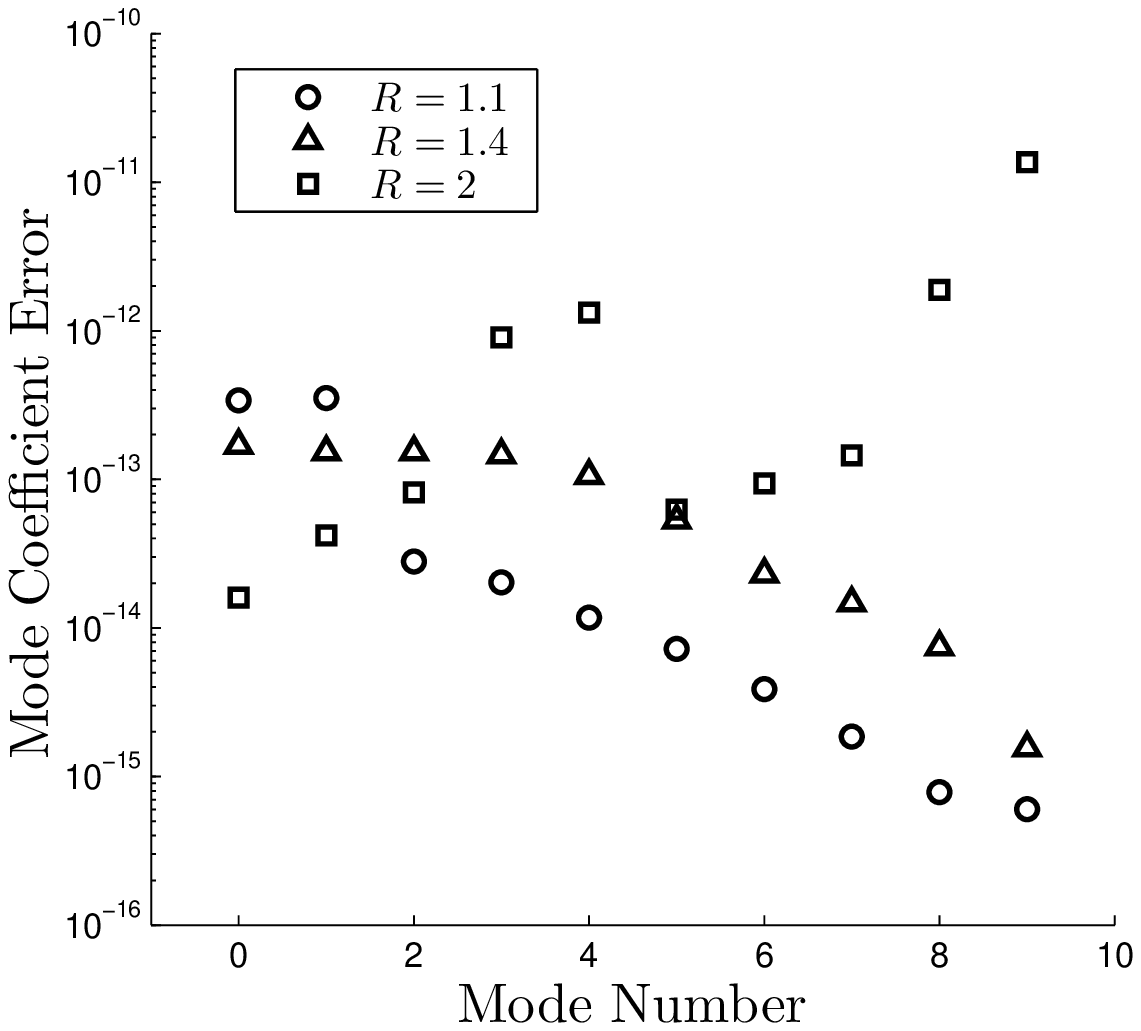}}
\caption{Maximum error in mode coefficients.}\label{fig:free_stream_b_err}
 \end{minipage}
 \hspace{0.5cm}
 \begin{minipage}[t]{0.5\linewidth}
\centerline{\includegraphics[height=6.1cm]{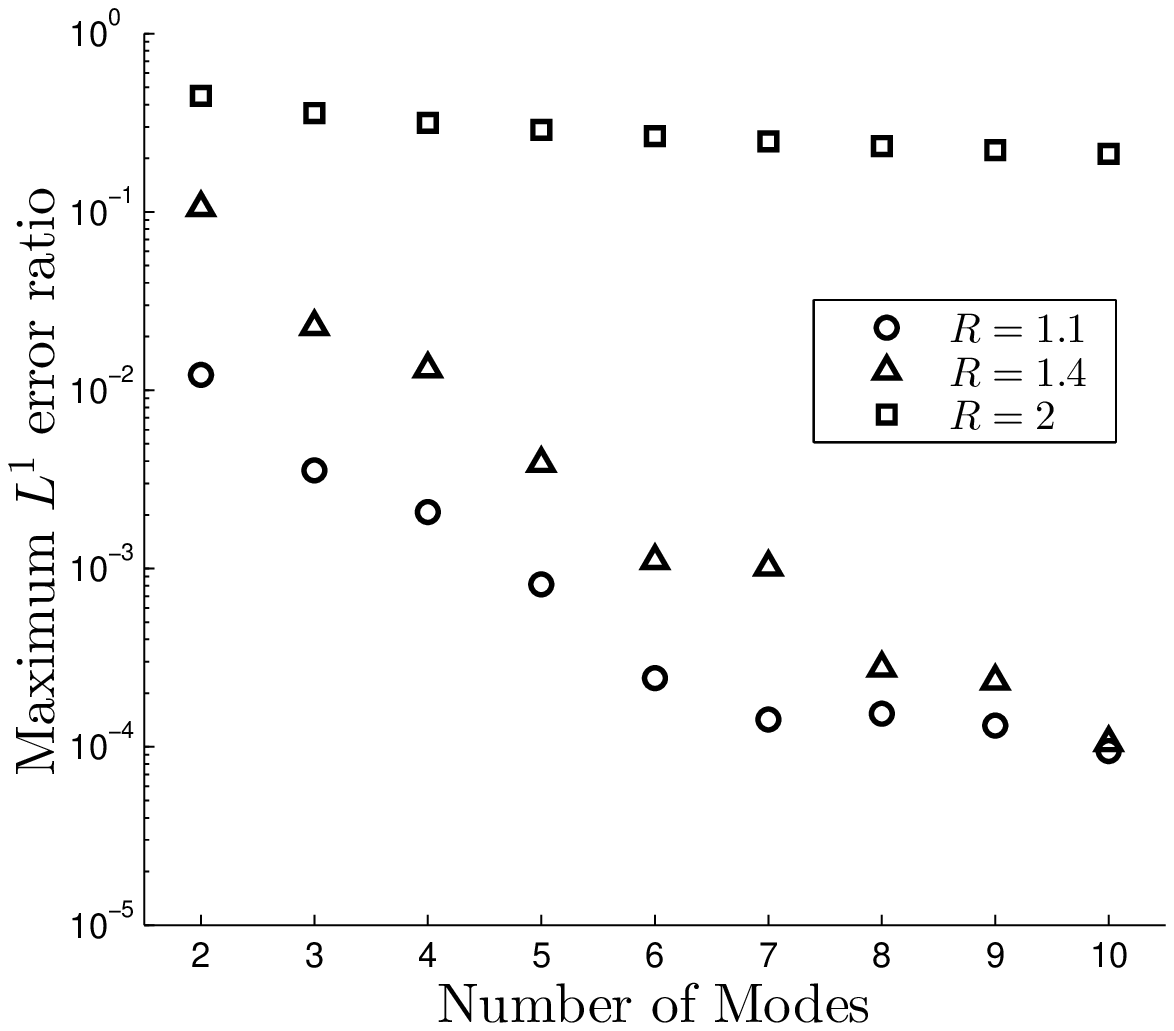}}
\caption{Maximum ratio  of $L^1$ error between computed and exact solutions to $L^1$ norm of the exact solution.}\label{fig:free_stream_L1_err}
 \end{minipage}
 \end{figure}

 To show that the numerical integration accurately captures the mode coefficients of the exact solution, \req{exact_sol}, we show the error between the computed coefficients and actual coefficients, denoted by $\tilde b_n$ and $b_n$ respectively
\begin{equation}\label{mode_err_def}
\text{error}_n=\max_{t} |\tilde{b}_n(t)-b_n(t)|,
\end{equation}
 in figure \ref{fig:free_stream_b_err}, where the evolution of the system was computed using $N=10$ modes.

\begin{figure}[H]
\begin{minipage}[t]{0.5\linewidth}
\centerline{\includegraphics[height=6.1cm]{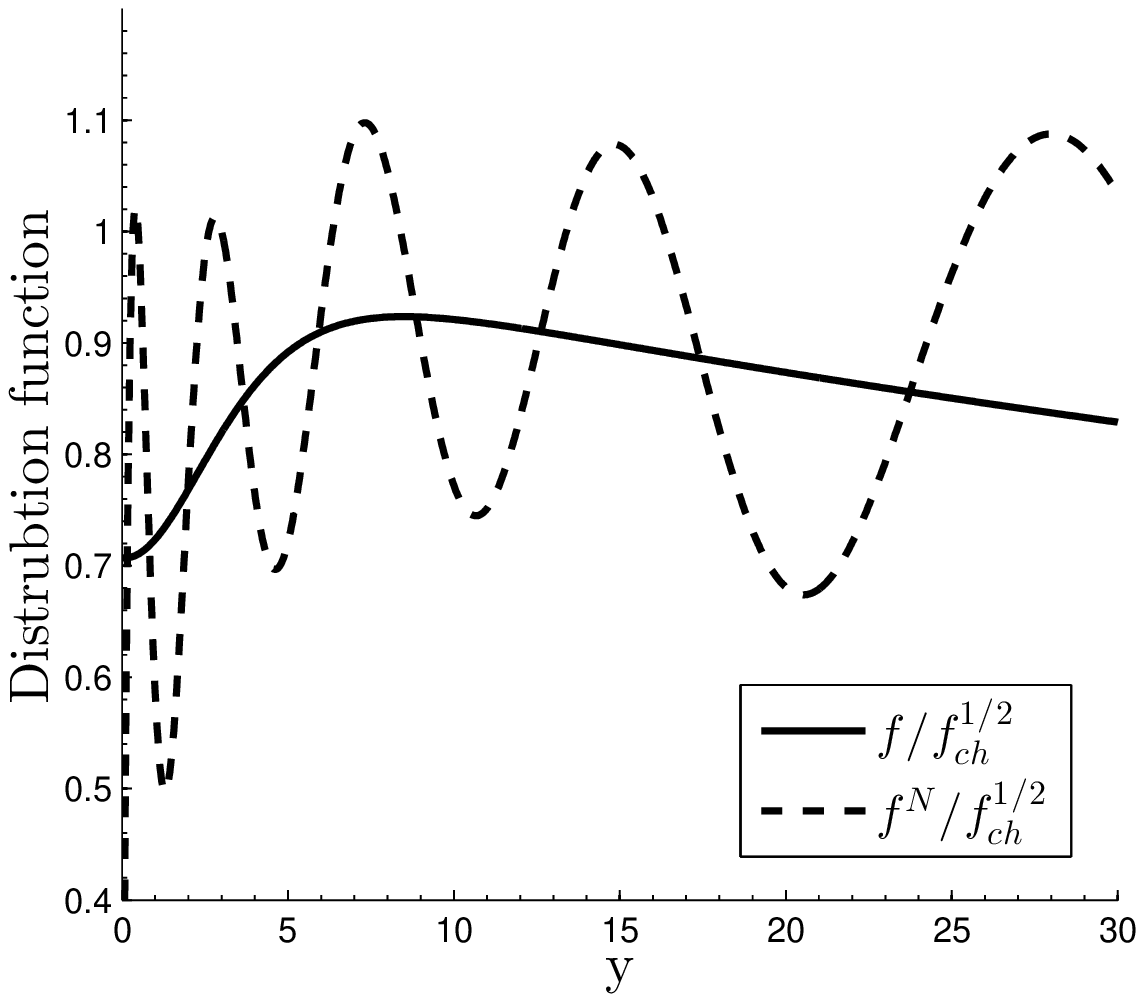}}
\caption{Approximate and exact solution for a reheating ratio $R=2$ and $N=10$ modes.}\label{fig:free_stream_approx_T_r_2}
 \end{minipage}
 \hspace{0.5cm}
 \begin{minipage}[t]{0.5\linewidth}
\centerline{\includegraphics[height=6.1cm]{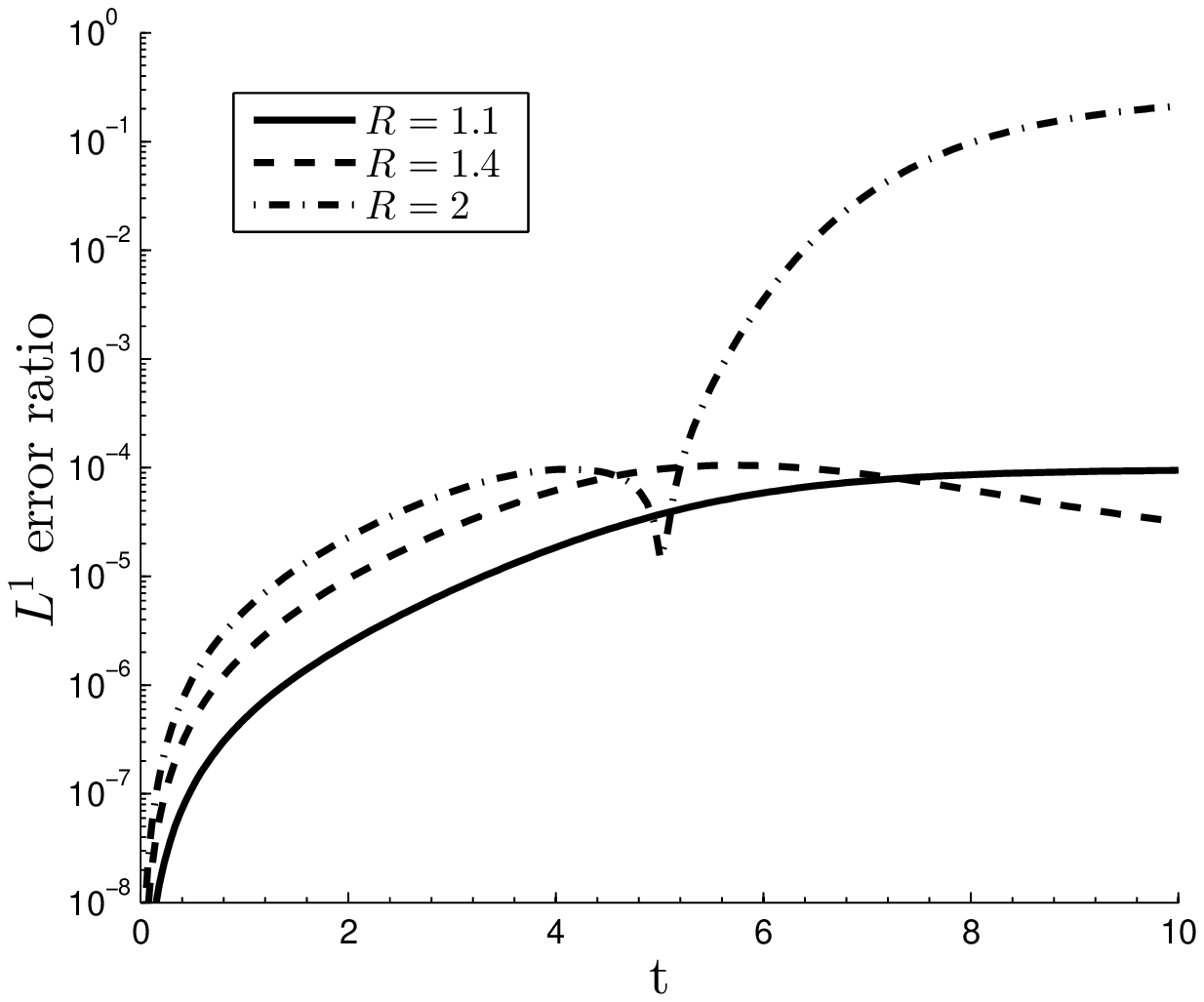}}
\caption{$L^1$ error ratio as a function of time for $N=10$ modes.}\label{fig:free_stream_L1_err_time}
 \end{minipage}
 \end{figure}

In figure  \ref{fig:free_stream_L1_err} we show the error between the exact solution $f$, and the numerical solution $f^N$ computed using $N=2,...,10$ modes over the solution time interval, where we define the error by
\begin{equation}\label{f_err}
\text{error}_N=\max_{t} \frac{\int |f-f^N|dy}{\int |f|dy}.
\end{equation}
For $R=1$ and $R=1.4$  the chemical equilibrium method works reasonably well (as long as the number of modes is at least 4, so that the energy and number densities are properly captured) but for $R=2$ the approximate solution exhibits spurious oscillations, as seen in figure \ref{fig:free_stream_approx_T_r_2}, and has poor $L^1$ error.  This is expected based on the prior results obtained in section \ref{basis_comparison}, as the exact solution is not in the Hilbert space used by this method for $R\geq 2$.  This is even clearer in figure \ref{fig:free_stream_L1_err_time} where we show the $L^1$ error ratio as a function of time for $N=10$ modes. In the $R=2$ case we see that the error increases as the reheating ratio approaches its asymptotic value of $R=2$ as $t\rightarrow\infty$.  As we will see, our methods achieves a much higher accuracy for a small number of terms in the case of large reheating ratio due to the replacement of dilution temperature scaling with the dynamical effective temperature $T$.  

\subsubsection{Chemical Non-Equilibrium Method}
We now solve  \req{free_stream_toy} using the chemical non-equilibrium method, with the orthonormal basis defined by the weight function \req{weight} for $N=2,...,10$ modes, a prescribed numerical integration tolerance of $10^{-13}$, and asymptotic reheating ratios of $R=1.1$, $R=1.4$, and $R=2$.  Recall that we are referring to $T$ and $\Upsilon$ as the first two modes ($n=0$ and $n=1$).

\begin{figure}[H]
 \begin{minipage}[b]{0.5\linewidth}
\centerline{\includegraphics[height=6.2cm]{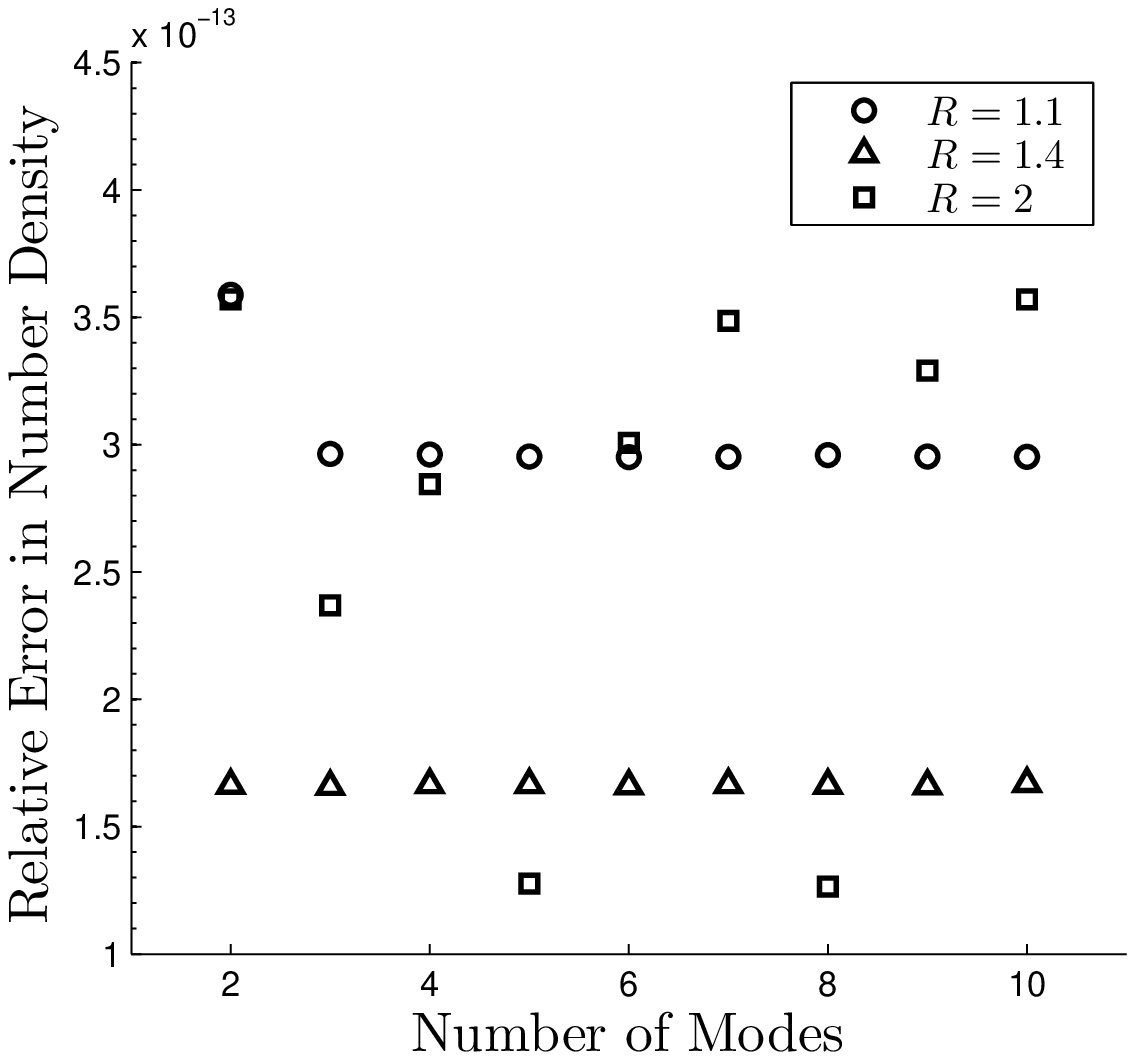}}
\caption{Maximum relative error in particle number density.}\label{fig:keq_num_err}
 \end{minipage}
 \hspace{0.5cm}
 \begin{minipage}[b]{0.5\linewidth}
\centerline{\includegraphics[height=6.2cm]{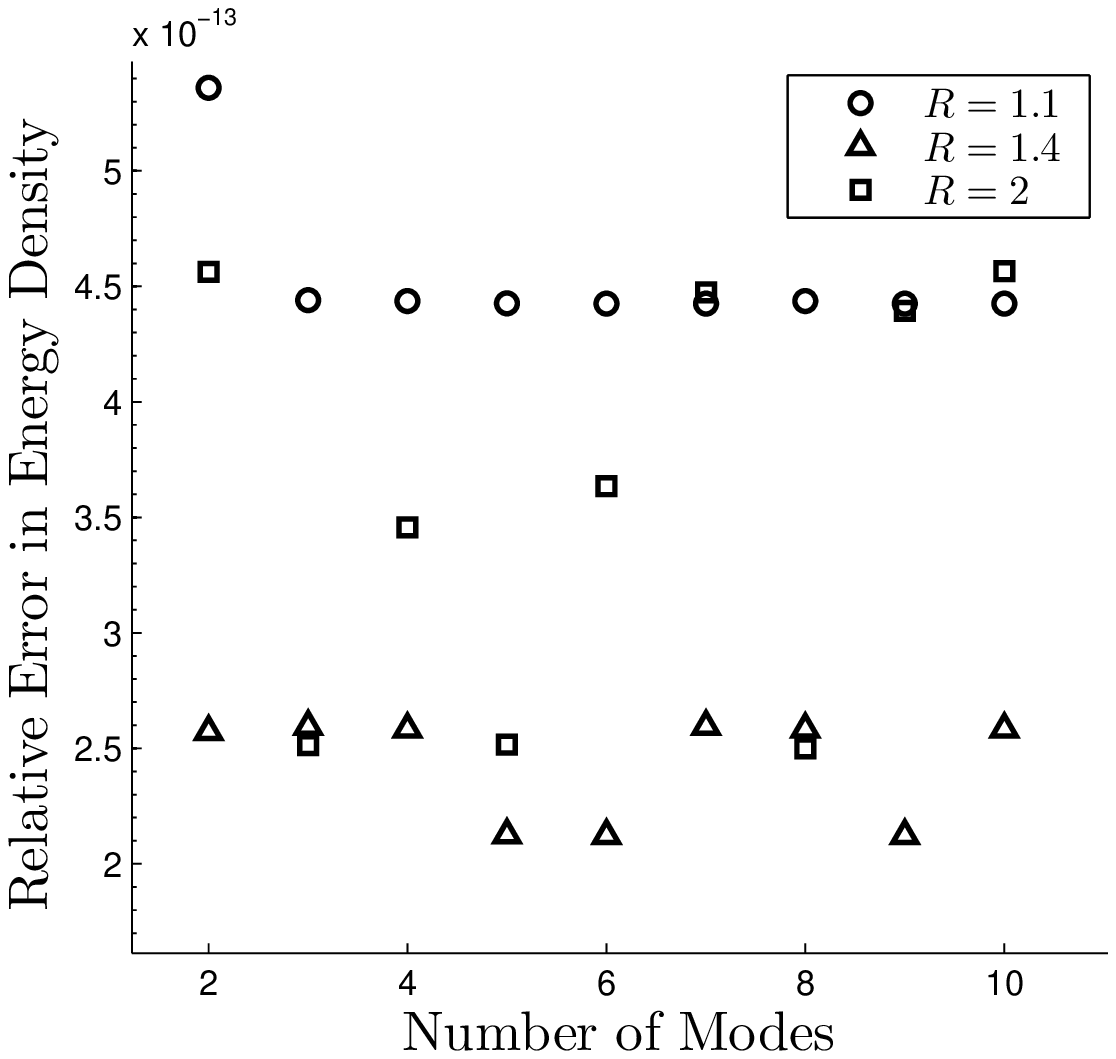}}
\caption{Maximum relative error in energy density.}\label{fig:keq_E_err}
 \end{minipage}
 \end{figure}

In figures \ref{fig:keq_num_err} and  \ref{fig:keq_E_err} we show the maximum relative error over the time interval $[0,10]$ in the number densities and energy densities respectively for various numbers of computed modes. Even for only $2$ modes, the number and energy densities are accurate up to the integration tolerance level.  This is in agreement with the analytical expressions in \req{th_eq_moments}.

\begin{figure}[H]
\begin{minipage}[t]{0.5\linewidth}
\centerline{\includegraphics[height=6cm]{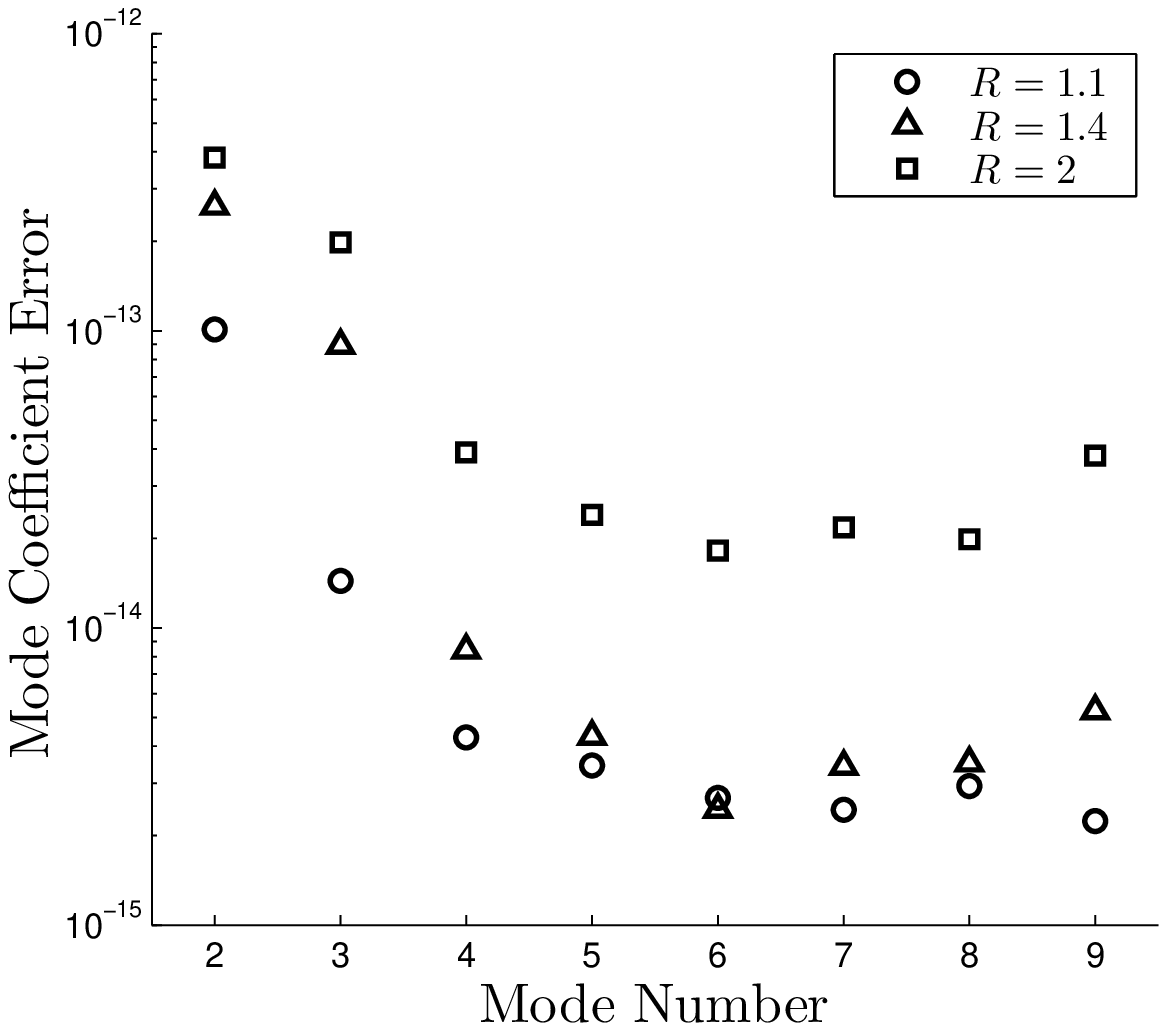}}
\caption{Maximum error in mode coefficients.}\label{fig:keq_b_err}
 \end{minipage}
 \hspace{0.5cm}
 \begin{minipage}[t]{0.5\linewidth}
\centerline{\includegraphics[height=6.1cm]{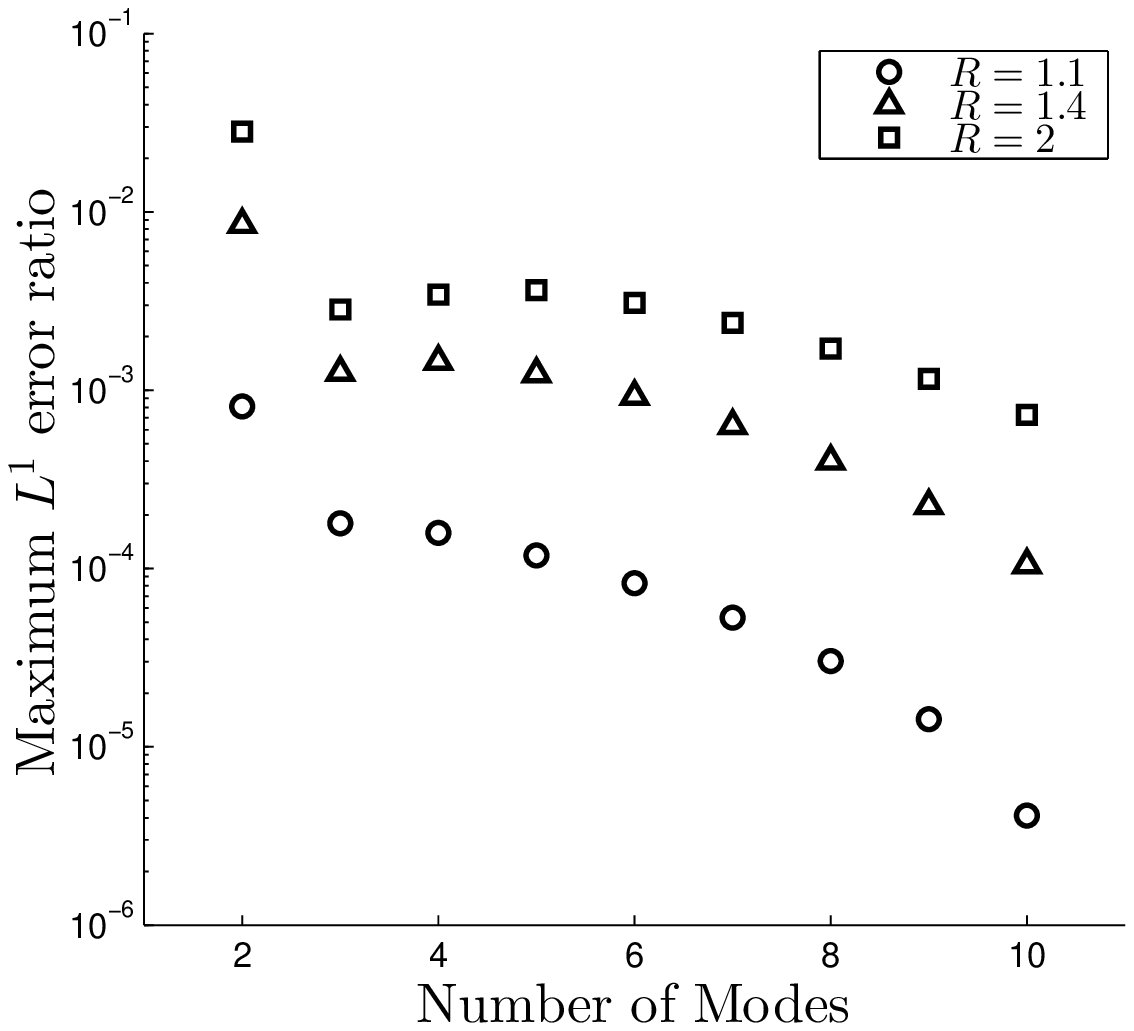}}
\caption{Maximum ratio of $L^1$ error between computed and exact solutions to $L^1$ norm of the exact solution.}\label{fig:keq_L1_err}
 \end{minipage}
 \end{figure}

 To show that the numerical integration accurately captures the mode coefficients of the exact solution, \req{exact_sol}, we give the error in the computed mode coefficients \req{mode_err_def}, where the evolution of the system was computed using $N=10$ modes, in figure \ref{fig:keq_b_err}. 

\begin{figure}[H]
 \begin{minipage}[t]{0.5\linewidth}
\centerline{\includegraphics[height=6.2cm]{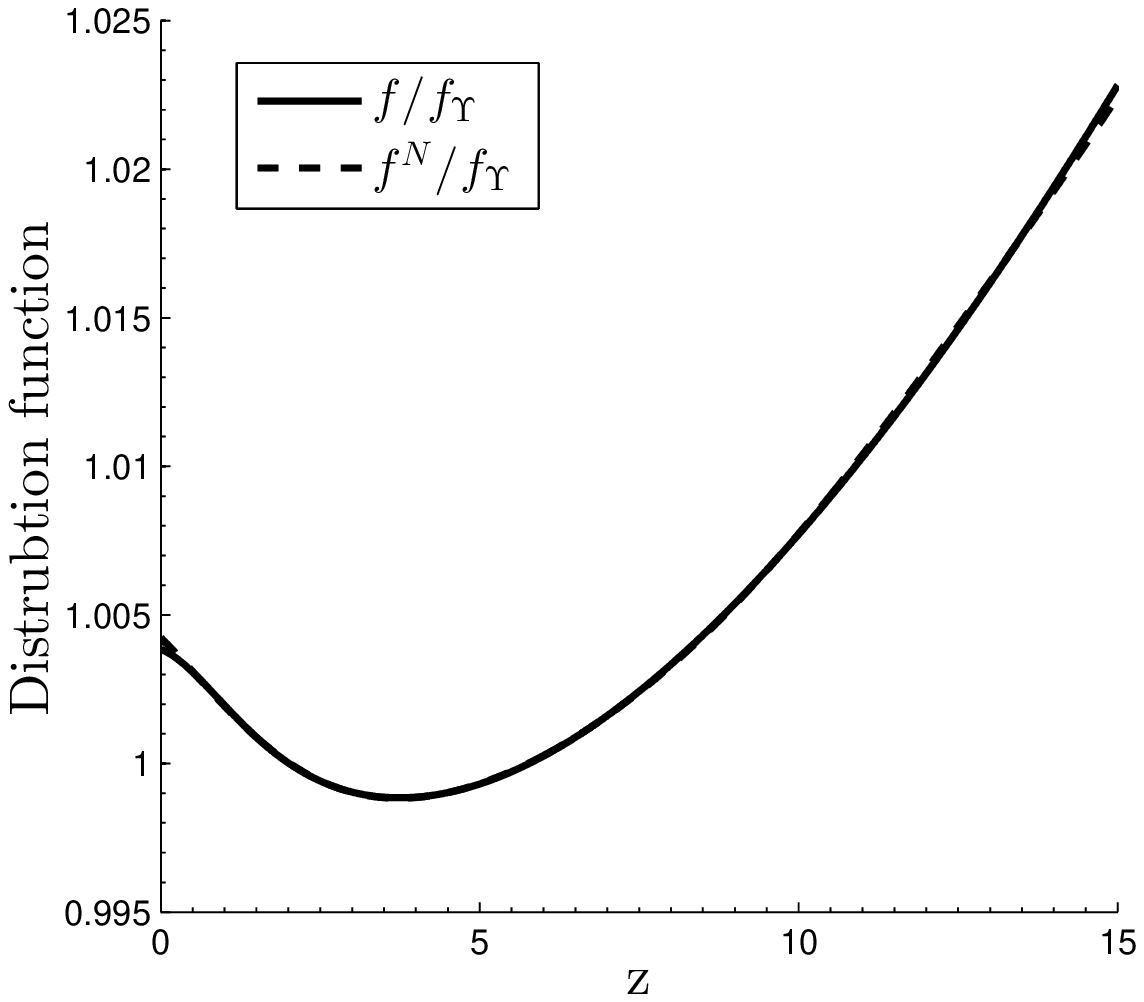}}
\caption{Approximate and exact solution for $R=2$ obtained with two modes.}\label{fig:keq_approx_Tr_2}
 \end{minipage}
 \hspace{0.5cm}
 \begin{minipage}[t]{0.5\linewidth}
\centerline{\includegraphics[height=6.2cm]{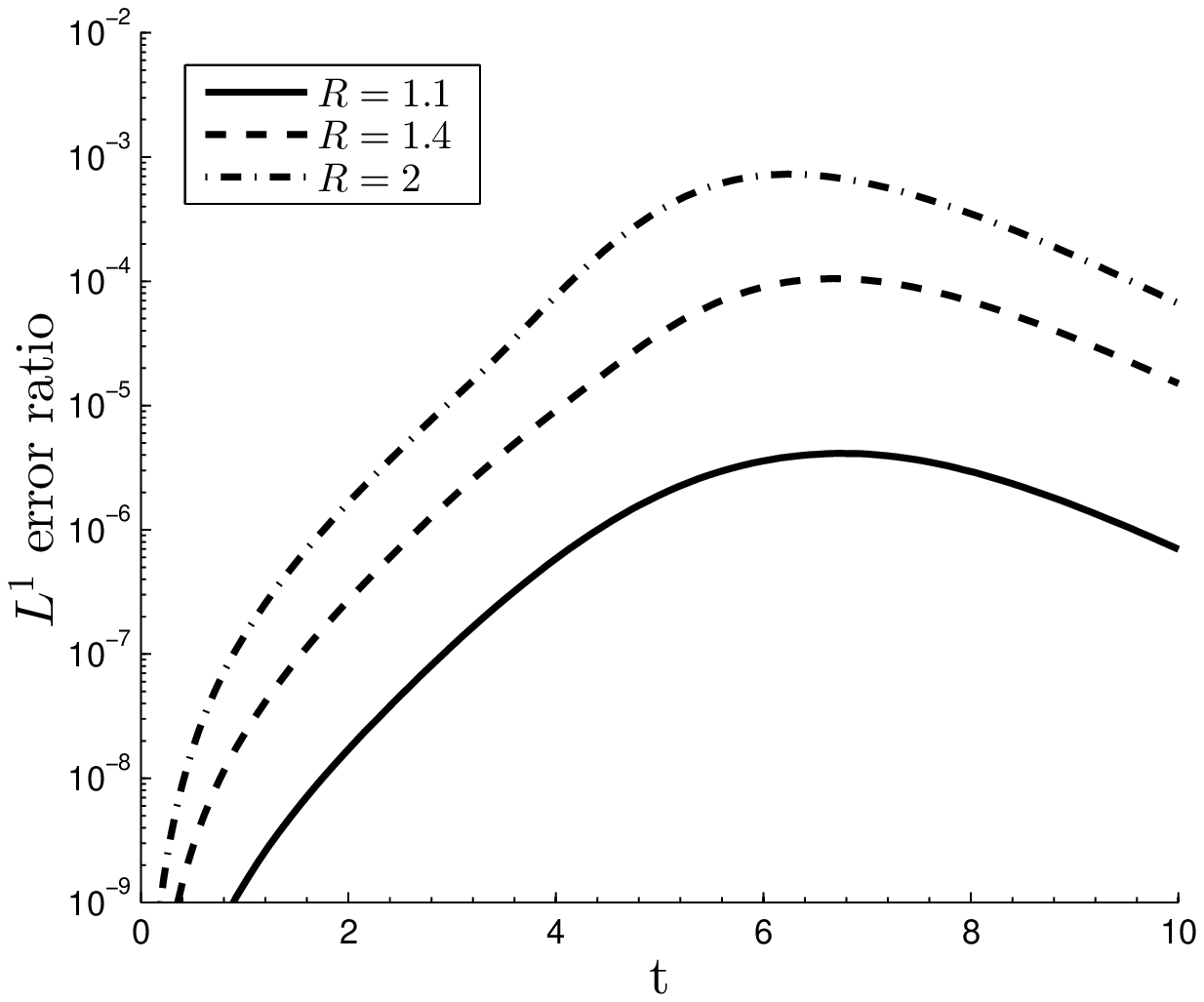}}
\caption{$L^1$ error ratio as a function of time for $n=10$ modes.}\label{fig:keq_L1_err_time}
 \end{minipage}
 \end{figure}

In figure \ref{fig:keq_L1_err} we show the error between the approximate and exact solutions, computed as in \req{f_err} for $N=2,...,10$ and $R=1.1$, $R=1.4$, and $R=2$ respectively.  For most mode numbers and $R$ values, the error using $2$ modes is substantially less than the error from the chemical equilibrium method using $4$ modes.  The result is most dramatic for the case of large reheating, $R=2$, where the spurious oscillations from the chemical equilibrium solution are absent, as seen in figure \ref{fig:keq_approx_Tr_2}, as compared to the chemical equilibrium method in figure \ref{fig:free_stream_approx_T_r_2}.  Note that we plot from $z\in [0,15]$ in comparison to $y\in[0,30]$ in figure \ref{fig:keq_approx_Tr_2} due to the relation $z=y/R$ as discussed in section \ref{basis_comparison}. Additionally, the error no longer increases as $t\rightarrow\infty$, as it did for the chemical equilibrium method, see figure \ref{fig:keq_L1_err_time}.  In fact it decreases since the exact solution approaches chemical equilibrium at a reheated temperature and hence can be better approximated by $f_\Upsilon$. 

In summary, in addition to the reduction in the computational cost when going from $4$ to $2$ modes, we also reduce the error compared to the chemical equilibrium method, all while still capturing the number and energy densities.  We emphasize that the error in the number and energy densities is limited by the integration tolerance and {\emph not} the number of modes (so long as the neglected, higher modes have a negligible impact on the first two modes, as is often the case).  

\subsubsection{Chemical Non-equilibrium Attractor}
The model problem in the previous subsections with $\Upsilon=1$  is reminiscent of the way that chemical equilibrium can emerge from chemical non-equilibrium in practice; the distribution of interest is attracted to some chemical equilibrium distribution but the particle creation/annihilation processes are not able to keep up with the momentum exchange and maintain an equilibrium particle yield, and so a fugacity $\Upsilon<1$ develops.  However, in order to isolate the effects of the fugacity on the solutions, we will now solve \req{toy_eq} under the  condition where our distribution is attracted to a fixed chemical non-equilibrium distribution.  More specifically, we take $T_{eq}(t)=1/a(t)$ and fugacities $\Upsilon=1.5$, $\Upsilon=0.9$, $\Upsilon=0.75$, and $\Upsilon=0.5$ with $a(t)$ defined as in \req{a_T_def}.

\begin{figure}[H]
 \begin{minipage}[t]{0.5\linewidth}
\centerline{\includegraphics[height=6.1cm]{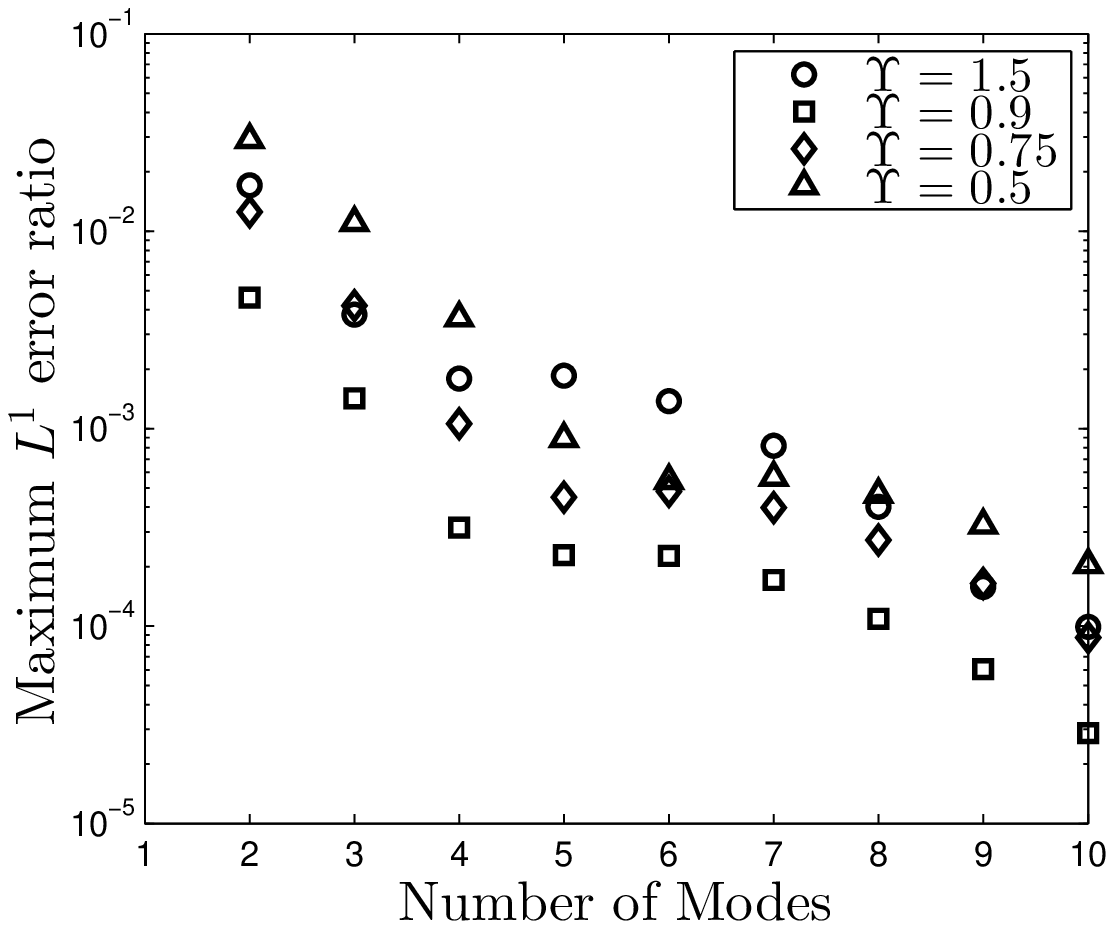}}
\caption{Maximum ratio  of $L^1$ error between computed and exact solutions to $L^1$ norm of the exact solution using the chemical equilibrium method.}\label{fig:free_stream_L1_err_Ups}
 \end{minipage}
 \hspace{0.5cm}
 \begin{minipage}[t]{0.5\linewidth}
\centerline{\includegraphics[height=6.1cm]{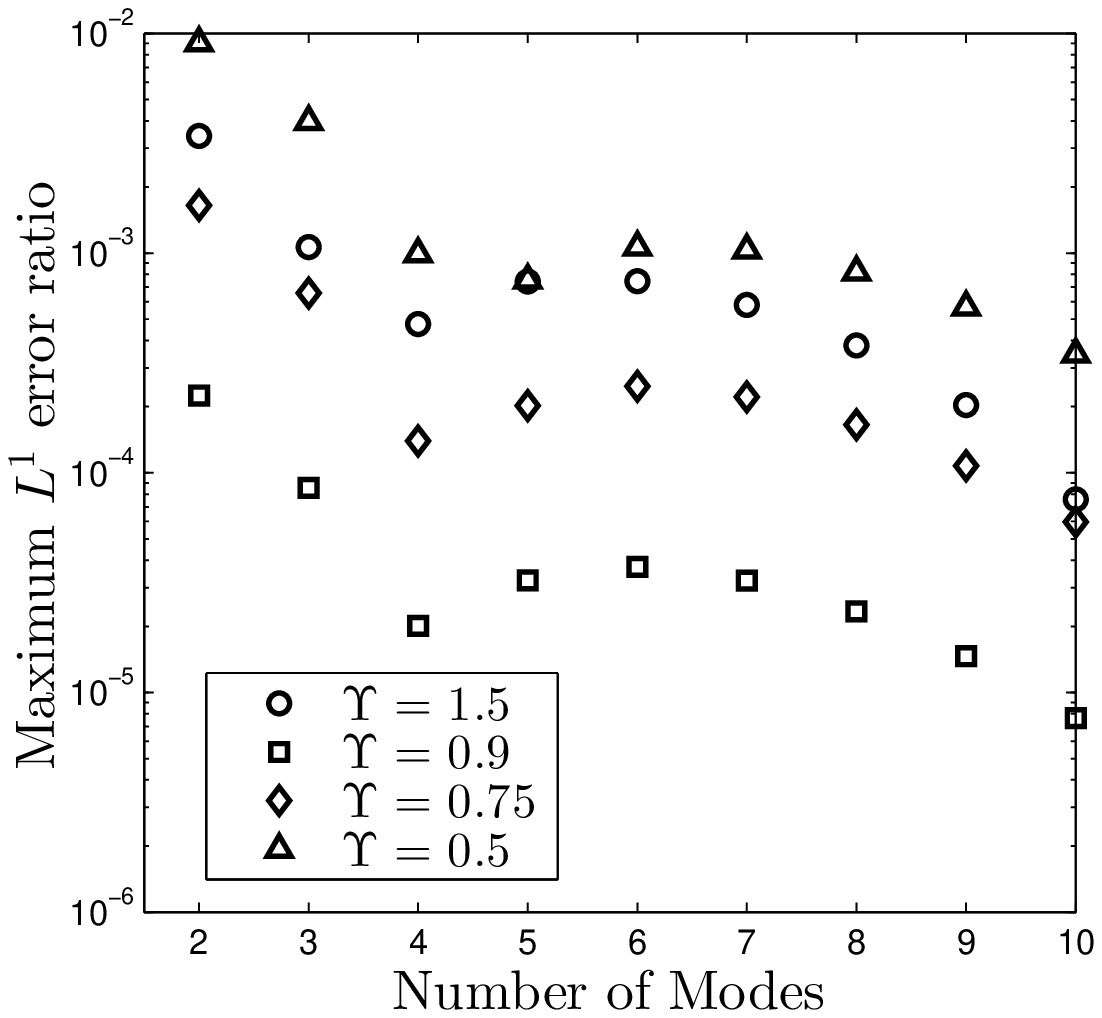}}
\caption{Maximum ratio of $L^1$ error between computed and exact solutions to $L^1$ norm of the exact solution using the chemical non-equilibrium method.}\label{fig:k_eq_L1_err_Ups}
 \end{minipage}
 \end{figure}

The behavior of the energy and particle number density errors are essentially the same as in the reheating test presented above and the mode coefficients are accurately captured, so we do not show these quantities here.  Instead  we show the maximum $L^1$ error, computed as in \req{f_err} in figure \ref{fig:free_stream_L1_err_Ups} for the chemical equilibrium method and in figure \ref{fig:k_eq_L1_err_Ups} for the chemical non-equilibrium method.  The error when using the latter method with only two terms is comparable to the former with four terms.

\begin{figure}[H]
\centerline{\includegraphics[height=6.1cm]{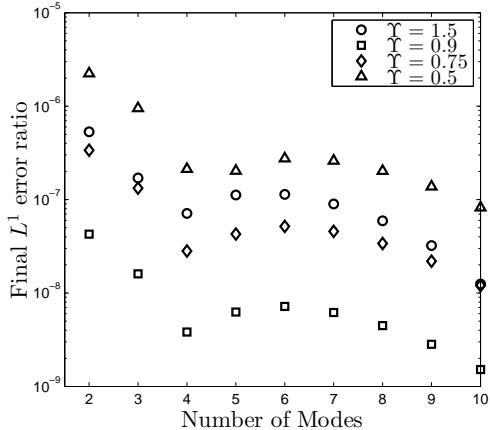}}
\caption{Final value of the ratio  of $L^1$ error between computed and exact solutions to $L^1$ norm of the exact solution using the chemical non-equilibrium method.}\label{fig:k_eq_L1_err_Ups_final}
 \end{figure}

In figure \ref{fig:k_eq_L1_err_Ups_final} we show the final value of the relative $L^1$ error,
\begin{equation}\label{f_err_final}
\text{error}_{N}= \frac{\int |f-f^N|dy}{\int |f|dy}.
\end{equation}
where the distribution functions are evaluated at $t=t_f=10$. This figure is for the chemical non-equilibrium method applied to the chemical non-equilibrium attractor. The error tolerance for numerical integration was chosen small enough that the error is determined by the number of modes.  The chemical non-equilibrium ansatz is able to represent the final asymptotic state accurately and so the final error is much smaller than the maximum error in figure \ref{fig:k_eq_L1_err_Ups}.  For the chemical equilibrium ansatz the final error is nearly the same as the maximum error and so we don't show it here.

\section{Summary and Outlook}
We have presented a spectral method for solving the Relativistic Boltzmann equation for a system  of massless fermions  diluting in time based on a dynamical basis of orthogonal polynomials.  The method is adapted to systems evolving near kinetic equilibrium, but allows for potentially strong chemical non-equilibrium in a transient and/or final state as well as strong reheating i.e. decoupling of temperature scaling from dilution scaling.  

The method depends on two time dependent parameters, the effective temperature $T(t)$ and phase space occupancy or fugacity $\Upsilon(t)$, whose dynamics are isolated by the requirement that the lowest modes capture the energy and particle number densities.  This gives the method a natural physical interpretation.  In particular, the dynamical fugacity is capable of naturally expressing the emergence of chemical non-equilibrium during the freeze-out process while the effective temperature captures any reheating phenomenon. Any system in approximate kinetic equilibrium that undergoes reheating and/or transitions to chemical non-equilibrium is a good match for this method.  In fact it is almost assured that our method will be considerably more computationally economical than the chemical equilibrium spectral method for any physical system in which the cost of computing the collision terms is high.

We validated the method on a model problem that exhibits the physical characteristics of reheating and chemical non-equilibrium.  We demonstrated that particle number and energy densities are captured accurately using only two degrees of freedom, the effective temperature and fugacity.  In general, this will hold so long as the back reaction from non-thermal distortions is small i.e. as long as kinetic equilibrium is a good approximation.

The method presented here should be compared to the spectral method used in \cite{Esposito2000,Mangano2002}, which uses a fixed basis of orthogonal polynomials and is adapted to systems that are close to chemical equilibrium (or with a non-dynamical chemical potential in \cite{Esposito2000}) with dilution temperature scaling.  In addition to more closely mirroring the physics of systems that exhibit reheating and chemical non-equilibrium, the method presented here has a computational advantage over the chemical equilibrium method. Even when the system is close to chemical equilibrium with dilution temperature scaling, as is the case for the problem studied in \cite{Esposito2000,Mangano2002}, the method presented here reduces the minimum number of degrees of freedom needed to capture the particle number and energy densities from four to two.  In turn, this reduces the minimum number of collision integrals that must be evaluated by more than half.  

Numerical evaluation of collision operators for realistic interactions is a costly operation and so the new `emergent chemical non-equilibrium' approach we have presented here constitutes a significant reduction in the numerical cost of obtaining solutions. Moreover, even if the chemical equilibrium approach were to be properly modified to gain mathematical advantages we show  in our chemical non-equilibrium approach, it is not at all clear that the chemical equilibrium method can, with comparable numerical effort, achieve  a precise solution under conditions where transient or final chemical non-equilibrium and reheating are strong.

Looking to future applications, the  gain in numerical efficiency we achieved should allow both space and time evolution to be considered in non-trivial dynamical models such as evolution of quark-gluon plasma fireball formed in relativistic heavy ion collisions. We expect to be able to explore  within the realm of Boltzmann dynamics    the question of `ideal' quark flow occurring at minimum viscosity~\cite{Romatschke:2007} and the shape of momentum distributions which systematically deviate from a thermal distribution~\cite{Wilk:2009}.  For study of the ensuing hadron flow  it is of relevance that we have been able to find suitable weights defining a spectrum of basis states capable of addressing the case of particles where the mass scale is relevant and chemical non-equilibrium is strong while retaining some of the advantages presented here for massless particles.  We will discuss the extension to the massive case in a future work.

\section*{Acknowledgment}
JB would like to acknowledge Government support under and awarded by DoD, Air Force Office of Scientific Research, National Defense Science and Engineering Graduate (NDSEG) Fellowship, 32 CFR 168a.  JW was supported in part by the US Department of Energy, Office of Science, Applied Scientific Computing Research, under award number DE-AC02-05CH11231, and by the National Science Foundation under award number DMS-0955078. JR was supported by a grant from the U.S. Department of Energy, DE-FG02-04ER41318.

\appendix

\section{Orthogonal Polynomials}\label{orthopoly_app}
\subsection{Generalities}\label{ortho-general}
Let $w:(a,b)\rightarrow [0,\infty)$ be a weight function where $(a,b)$ is a (possibly unbounded) interval and consider the Hilbert space $L^2(w dx)$.   We will consider weights such that $x^n\in L^2(wdx)$ for all $n\in\mathbb{N}$. We denote the inner product by $\langle\cdot,\cdot\rangle$, the norm by $||\cdot||$, and for a vector $\psi\in L^2$ we let $\hat{\psi}\equiv \psi/||\psi||$.  The classical three term recurrence formula can be used to define a set of orthonormal polynomials $\hat{\psi}_i$ using this weight function, for example see \cite{Olver},
\begin{align}\label{poly_recursion}
&\psi_0=1, \hspace{2mm} \psi_1=||\psi_0||(x-\langle x\hat\psi_0,\hat{\psi}_0\rangle)\hat{\psi}_0,\\
&\psi_{n+1}=||\psi_n||\left[\left(x-\langle x\hat\psi_n,\hat\psi_n\rangle\right)\hat\psi_n-\langle x\hat\psi_n,\hat\psi_{n-1}\rangle\hat\psi_{n-1}\right].
\end{align}
 In an effort to reduce the operation count for evaluation of the $\hat\psi_i$'s, one could use the above recurrence relation to generate the coefficients of the $\hat \psi_i$'s in the monomial basis and use evaluation schemes based these coefficients, however this can lead to accuracy issues.  For a comparison of the recursive formula to Horner's scheme for evaluating orthogonal polynomials see, for example, \cite{Fateman} which considers Legendre polynomials.  In the sections of this work where accuracy is paramount we have used the recursive formula given above.

 One can also derive recursion relations for the derivatives of $\psi_n$ with respect to $x$, denoted with a prime,
\begin{align}\label{deriv_rec}
&\psi_0^{'}=0, \hspace{2mm} \hat{\psi}_1^{'}=\frac{||\psi_0||}{||\psi_1||}\hat\psi_0,\\
&\hat{\psi}_{n+1}^{'}=\frac{||\psi_n||}{||\psi_{n+1}||}\left[\hat\psi_n+\left(x-\langle x\hat\psi_n,\hat\psi_n\rangle\right)\hat{\psi}_n^{'}-\langle x\hat\psi_n,\hat\psi_{n-1}\rangle\hat{\psi}_{n-1}^{'}\right].
\end{align}
Since $\hat{\psi}_n^{'}$ is a degree $n-1$ polynomial, we have the expansion 
\begin{equation}\label{app:psi_prime_exp}
\hat{\psi}_n^{'}=\sum_{k<n} a_n^k \hat{\psi}_k.
\end{equation}
Using \req{deriv_rec} we obtain a recursion relation for the $a_n^k$
\begin{equation}\label{app:psi_prime_coeff}
a_{n+1}^k=\frac{||\psi_n||}{||\psi_{n+1}||}\left(\delta_{n,k}-\langle x\hat\psi_n,\hat\psi_{n}\rangle a_n^k-\langle x\hat\psi_n,\hat\psi_{n-1}\rangle a_{n-1}^k+\sum_{l<n}a_n^l\langle x\hat\psi_l,\hat\psi_k\rangle\right),\notag
\end{equation}
\begin{equation}
a_1^0=\frac{||\psi_0||}{||\psi_1||}.
\end{equation}
This can be used to compute some of the inner products appearing in \req{A_B_matrices}.
\subsection{Parametrized Families of Orthogonal Polynomials}\label{ortho-polynom-fam}
Our method requires not just a single set of orthogonal polynomials, but rather a parametrized family of orthogonal polynomial generated by a weight function $w_t(x)$ that is a $C^1$ function of both $x\in(a,b)$ and some parameter $t$.  To emphasize this, we write $g_t(\cdot,\cdot)$ for $\langle\cdot,\cdot\rangle$.  We will assume that $\partial_t w$ is dominated by some $L^1(dx)$ function of $x$ only that decays exponentially as $x\rightarrow\pm\infty$ (if the interval is unbounded). In particular, this holds for the weight function \req{weight}.

Given the above assumption about the decay of $\partial_t w$, the dominated convergence theorem implies that $\langle p,q\rangle$ is a $C^1$ function of $t$ for all polynomials $p$ and $q$ and justifies  differentiation under the integral sign. By induction, it also implies implies that the $\hat\psi_i$ have coefficients that are $C^1$ functions of $t$. Therefore, for any polynomials $p$, $q$ whose coefficients are $C^1$ functions of $t$ we have
\begin{equation}
\frac{d}{dt}g_t( p,q)=\dot{g}_t(p,q)+g_t(\dot{p},q)+g_t( p,\dot{q})
\end{equation}
where a dot denotes differentiation with respect to $t$ and we use $\dot{g}_t(\cdot,\cdot)$ to denote the inner product with respect to the weight $\dot{w}$.  

Ostensibly, \req{b_eq} for the mode coefficients requires us to compute $g(\dot{\hat\psi}_i,\hat\psi_j)$ for all $i,j$.  However, as we will show in the following section, only the result for $i=j$ is needed.  This can be obtained quite easily by differentiating the relation
\begin{equation}
\delta_{ij}=g_t(\hat\psi_i,\hat\psi_j),
\end{equation}
yielding
\begin{equation}\label{ortho_deriv_eq}
0=\dot g_t(\hat\psi_i,\hat\psi_j)+g_t(\dot{\hat\psi}_i,\hat\psi_j)+g_t(\hat\psi_i,\dot{\hat\psi}_j).
\end{equation}
For $i=j$ we obtain
\begin{equation}\label{norm_deriv_eq}
g_t(\dot{\hat\psi}_i,\hat\psi_i)=-\frac{1}{2}\dot{g}_t(\hat\psi_i,\hat\psi_i).
\end{equation}

Though they will not be needed, we give the results for $i\neq j$ for completeness.  Since $\hat\psi_j$ is orthogonal to all polynomials of degree less than $j$, for $j>i$ 
\begin{equation}
g_t(\dot{\hat\psi}_i,\psi_j)=0
\end{equation}
and for $j<i$ we have
\begin{equation}
g_t(\dot{\hat\psi}_i,\psi_j)=-\dot{g}_t(\hat\psi_i,\hat\psi_j)-g(\hat\psi_i,\dot{\hat\psi}_j)=-\dot{g}_t(\hat\psi_i,\hat\psi_j).
\end{equation}

\subsection{Proof of Lower Triangularity}\label{lower_triang}
Here we prove that the matrices that define the dynamics of the mode coefficients $b^k$ are lower triangular.  Recall the definitions
\begin{align}
A^k_i(\Upsilon)\equiv&\langle\frac{z}{f_\Upsilon }\hat\psi_i\partial_zf_\Upsilon ,\hat\psi_k\rangle+\langle z\partial_z \hat\psi_i,\hat\psi_k\rangle,\label{app:A_def}\\
B^k_i(\Upsilon)\equiv &\Upsilon\left(\langle\frac{1}{f_\Upsilon }\frac{\partial f_\Upsilon }{\partial\Upsilon}\hat\psi_i,\hat\psi_k\rangle+\langle\frac{\partial\hat{\psi}_i}{\partial \Upsilon},\hat\psi_k\rangle\right).\label{app:B_def}
\end{align}
Using integration by parts, we see that
\begin{equation}
A^k_i=-3\langle\hat\psi_i,\hat\psi_k\rangle-\langle \hat \psi_i,z\partial_z\hat\psi_k\rangle.
\end{equation}
Since $\hat\psi_i$ is orthogonal to all polynomials of degree less than $i$ we have $A^k_i=0$ for  $k<i$.  For $k>i$ the second term in \req{app:A_def} vanishes and so the only terms involving $\partial_z\hat\psi_i$ that are needed are the diagonal entries $\langle z\partial_z \hat\psi_k,\hat\psi_k\rangle$ which, using \req{app:psi_prime_exp} and \req{app:psi_prime_coeff}, can be simplified to
\begin{equation}
\langle z\partial_z \hat\psi_k,\hat\psi_k\rangle=\sum_{l<k}a_k^l\langle z\hat\psi_l,\hat\psi_k\rangle=a_k^{k-1}\langle z\hat\psi_{k-1},\hat\psi_k\rangle\delta_{k-1\geq 0}.
\end{equation}

$B^k_i$ can be simplified as follows.  First differentiate 
\begin{equation}
\delta_{ik}=\langle \hat\psi_i,\hat\psi_j\rangle
\end{equation}
with respect to $\Upsilon$ to obtain
\begin{align}
0=&\int \hat\psi_i\hat\psi_k\partial_{\Upsilon}wdz+\langle \partial_{\Upsilon}\hat\psi_i,\hat\psi_k\rangle+\langle \hat\psi_i,\partial_{\Upsilon}\hat\psi_k\rangle\\
=&\langle\frac{\hat\psi_i}{f_\Upsilon}\partial_{\Upsilon}f_\Upsilon,\hat\psi_k \rangle+\langle\partial_{\Upsilon}\hat\psi_i,\hat\psi_k\rangle+\langle \hat\psi_i,\partial_{\Upsilon}\hat\psi_k\rangle
\end{align}
Therefore 
\begin{equation}
B^k_i=-\Upsilon\langle\hat\psi_i,\partial_{\Upsilon}\hat\psi_k\rangle.
\end{equation}
$\partial_\Upsilon \hat\psi_k$ is a degree $k$ polynomial, hence $B_i^k=0$ for $k<i$ as claimed. For $k>i$  the second term in the definition of $B_i^k$, \req{app:B_def}, vanishes by the same reasoning.  Therefore, the only inner products involving $\partial_\Upsilon\hat{\psi}_i$ that are required are  
\begin{equation}
B_k^k=-\Upsilon\langle\hat\psi_k,\partial_{\Upsilon}\hat\psi_k\rangle.
\end{equation}
We showed how to compute these in \ref{ortho-polynom-fam}.


\end{document}